\documentclass[12pt,a4paper]{amsart}
\usepackage{a4wide}
\usepackage[english,russian]{babel}
\usepackage[T2A]{fontenc}
\usepackage[cp1251]{inputenc} 
\usepackage{amsfonts}
\usepackage{amssymb, amsthm, amscd}
\usepackage{amsmath}
\usepackage{mathtools}
\usepackage{needspace}
\usepackage{etoolbox}
\usepackage{lipsum}
\usepackage{comment}
\usepackage{cmap}
\usepackage[pdftex]{graphicx}
\usepackage[unicode]{hyperref}
\usepackage[matrix,arrow,curve]{xy}
\usepackage[usenames,dvipsnames]{xcolor}
\usepackage{colortbl}
\usepackage{textcomp}
\usepackage{cite}
\usepackage{euscript}

\DeclareMathOperator{\Lie}{Lie}

\DeclareMathOperator{\ring}{ring}
\DeclareMathOperator{\aug}{aug}
\DeclareMathOperator{\ad}{ad}
\DeclareMathOperator{\fppf}{fppf}

\newcommand{\Z}{\mathbb{Z}}

\newcommand{\N}{\mathbb{N}}

\newcommand{\Root}{\mathfrak{R}}

\renewcommand{\C}{\mathbb{C}}

\newcommand{\A}{\mathfrak{A}}
\newcommand{\B}{\mathfrak{B}}

\newcommand{\ovl}{\overline}
\newcommand{\Equal}{\Leftrightarrow}

\newcommand{\eps}{\varepsilon}
\newcommand{\ph}{\varphi}

\renewcommand{\Im}{\mathop{\mathrm{Im}}\nolimits}

\newcommand{\sub}{\subseteq}
\newcommand{\Hom}{\mathop{\mathrm{Hom}}\nolimits}
\newcommand{\Ker}{\mathop{\mathrm{Ker}}\nolimits}
\newcommand{\lev}{\mathop{\mathrm{lev}}\nolimits}

\newcommand{\rk}{\mathop{\mathrm{rk}}\nolimits}

\newcommand{\fact}[2]{#1/#2}

\renewcommand{\ge}{\geqslant}
\renewcommand{\le}{\leqslant}
\newcommand{\sm}{\setminus}
\newcommand{\map}[3]{#1\colon #2\to #3}
\newcommand{\phan}{\phantom{,}}
\newcommand{\gap}{\;\;}

\newcommand{\restr}[2]{\left. #1\right|_{#2}}
\newcommand{\<}{\langle}
\renewcommand{\>}{\rangle}
\newcommand{\+}{\oplus}
\newcommand{\emp}{\varnothing}

\renewcommand{\P}{\EuScript{P}}
\newcommand{\euQ}{\EuScript{Q}}
\newcommand{\euL}{\EuScript{L}}

\theoremstyle{plain}
\newtheorem{thm}{Theorem}
\newtheorem{lem}{Lemma}
\newtheorem{prop}{Proposition}

\theoremstyle{definition}

\newtheorem{cor}{Corollary}

\theoremstyle{remark}

\AtBeginEnvironment{thm}{\begin{samepage}}
\AtEndEnvironment{thm}{\end{samepage}}
\AtBeginEnvironment{lem}{\begin{samepage}}
\AtEndEnvironment{lem}{\end{samepage}}
\AtBeginEnvironment{st}{\begin{samepage}}
\AtEndEnvironment{st}{\end{samepage}}
\AtBeginEnvironment{crit}{\begin{samepage}}
\AtEndEnvironment{crit}{\end{samepage}}
\AtBeginEnvironment{ax}{\begin{samepage}}
\AtEndEnvironment{ax}{\end{samepage}}
\AtBeginEnvironment{defn}{\begin{samepage}}
\AtEndEnvironment{defn}{\end{samepage}}
\AtBeginEnvironment{cor}{\begin{samepage}}
\AtEndEnvironment{cor}{\end{samepage}}
\AtBeginEnvironment{note}{\begin{samepage}}
\AtEndEnvironment{note}{\end{samepage}}
\AtBeginEnvironment{prop}{\begin{samepage}}
\AtEndEnvironment{prop}{\end{samepage}}

\newcommand{\Stab}{\mathop{\mathrm{Stab}}\nolimits}
\newcommand{\SL}{\mathop{\mathrm{SL}}\nolimits}

\newcommand{\Tran}{\mathop{\mathrm{Tran}}\nolimits}

\newcommand{\tc}{\text{,}}
\newcommand{\tp}{\text{.}}
\renewcommand{\tilde}{\widetilde}

\DeclareMathOperator{\SC}{sc}
\DeclareMathOperator{\gen}{gen}

\makeatletter
\def\@settitle{\begin{center}%
    \baselineskip14\p@\relax
    \bfseries
    \@title
  \end{center}%
}

\def\@evenhead{\hfil\sc p. gvozdevsky\hfil}
\def\@oddhead{\hfil\sc overgroups of Levi subgroups I.\hfil}
\makeatother

\title{OVERGROUPS OF LEVI SUBGROUPS I. THE CASE OF ABELIAN UNIPOTENT RADICAL}

\author{P.~Gvozdevsky}
\date{}

\address{Chebyshev Laboratory,\\ St. Petersburg \\ State University, \\14th Line V.O., 29B, \\
	Saint Petersburg 199178 Russia}

\thanks{This publication is supported by Russian Science Foundation grant \textnumero\,17-11-01261}

\begin{document}
\selectlanguage{english}

\maketitle

\begin{abstract}
In the present paper we prove sandwich classification for the overgroups of the subsystem subgroup $E(\Delta,R)$ of the Chevalley group $G(\Phi,R)$ for the three types of pairs $(\Phi,\Delta)$ (the root system and its subsystem) such that the group $G(\Delta,R)$is (up to torus) a Levi subgroup of the parabolic subgroup with abelian unipotent radical. Namely we show that for any such an overgroup $H$ here exists a unique pair of ideals $\sigma$ of the ring $R$ such that $E(\Phi,\Delta,R,\sigma)\le H\le N_{G(\Phi,R)}(E(\Phi,\Delta,R,\sigma))$.  
\end{abstract}

\section{Introduction}

In the paper \cite{Aschbacher84}, dedicated to the {\it maximal subgroups project}, Michael Aschbacher introduced eight classes $C_1$-$C_8$ 8 of subgroups of finite simple classical groups. The groups from these classes are ''obvious'' maximal subgroups of a finite classical groups. To be precise, each subgroup from an Aschbacher class either is maximal itself or is contained in a maximal subgroup that in its turn either also belongs to an Aschbacher class or can be constructed by a certain explicit procedure.

Nikolai Vavilov defined five classes of “large” subgroups of the Chevalley groups (including exceptional ones) over arbitrary rings (see \cite{StepDiss} for details). Although these
subgroups are not maximal, he conjectured that they are sufficiently large for the corresponding overgroup lattice to admit a description. One of these classes is the class of
subsystem subgroups (the definition will be given in Subsection 2.1).

Overgroups of (elementary) subsystem subgroups of the general linear group were studied in the papers \cite{BVNnets},\cite{BV82GL},\cite{BV82},\cite{BV84},\cite{VavGLSemiloc}.In that case, subsystem subgroups are precisely block-diagonal subgroups. This results were generalized to orthogonal and symplectic groups under the assumption $2\in R^*$ in the thesis of Nikolai Vavilov (see also \cite{VavMIAN}, \cite{VavSplitOrt2001} and \cite{VavSimpSubsyst}). After that, this assumption was lifted in the thesis of Alexander Shchegolev \cite{SchegDiss}. In
that thesis, Shchegolev also solved the problem for unitary groups (see also \cite{SchegMainResults} and \cite{SchegSymplectic}).

The problem we discuss in the present paper is to describe the overgroups of subsystem
subgroups of exceptional groups (over a commutative ring). This problem was posed
in \cite{VavStepSurvey} (Problem 7). The first step of solving this problem was done in \cite{VSch}.The table from that paper contains the list of pairs $(\Phi,\Delta)$ for which such a description may be possible in principle, along with the number of ideals determining the level and along with certain links between these ideals.

By a standard description we mean {\it ''sandwich classification''}. Let $G$ be an abstract
group, and let $\euL$ be a certain lattice of its subgroups. The lattice $\euL$  admits sandwich classification if it is a disjoint union of ''sandwiches'':

\begin{align*}
&\euL=\bigsqcup_i L(F_i,N_i)\tc\\
&L(F_i,N_i)=\{H\colon F_i\le H\le N_i\}\tc
\end{align*}

where $i$ runs through some index set, and $F_i$ is a normal subgroup of $N_i$ for all $i$. To study such a lattice, it suffices to study the quotients $N_i/F_i$. In \cite{VSch} it was conjectured that the lattice of subgroups of a Chevalley group that contain a sufficiently large subsystem subgroup admits sandwich classification.

In the present paper, we prove a sandwich classification theorem for the embeddings $A_{l-1}\le D_l$ (where $l\ge 5$), $D_5\le E_6$, and $E_6\le E_7$. Our approach allows us to consider
all three cases simultaneously. In all these cases, the corresponding subsystem subgroup
is (up to a torus) a Levi subgroup $L_{\alpha}$, where $\alpha$ is a fundamental root such that its coefficient in the decomposition of a maximal root is equal to one. In that case, the level
is determined by two ideals. One can add to this list the embedding $A_{l-1}\le A_l$, for which
the main result follows, of course, from the paper \cite{BV84}. The fact that in all these cases the corresponding unipotent radical is Abelian simplifies our problem from a technical point of view.

Overgroup description for the embeding $A_{l-1}\le D_l$ is a special case of the results
obtained in \cite{VavSplitOrt} and \cite{SchegDiss}. Another two cases considered in the present paper are new.

We mention two previous results that are closely related to the present paper.

\begin{itemize}
\item Over a field (distinct from $F_3$ and of characteristic not equal to 2), overgroup
description in the cases considered in the present paper was obtained by Wang Dengyin in \cite{WangDengyin}. His proof involves Bruhat decomposition.

\item Description of subgroups of a maximal parabolic subgroup that contain the elementary Levi subgroup was obtained in the paper of Anastasia Stavrova \cite{StavNormStr}. Note that over a field such a description had been obtained before in the papers of Wang Dengyin and Li ShangZhi \cite{WangDengyinLiShangZhi} and Victoria Kozakevich and Anastasia Stavrova \cite{KazStavr}.
\end{itemize} 

\section*{Acknowledgment}

I am grateful to my teacher Nikolai Vavilov for setting the problem and for extremely
helpful suggestions.

\section{Basic notation}

\subsection{Root systems and Chevalley groups}

Let $\Phi$ be an irreducible root system, $\P$ a lattice that is intermediate between the root lattice $\euQ(\Phi)$ and the weight lattice $\P(\Phi)$, $R$ a commutative associative ring with unity, $G(\Phi,R)=G_\P(\Phi,R)$ a Chevalley group of type $\Phi$ over $R$, $T(\Phi,R)=T_\P(\Phi,R)$  a split maximal torus of $G_\P(\Phi,R)$. For every root $\alpha\in\Phi$ we denote by $X_\alpha=\{x_\alpha(\xi)\colon \xi\in R\}$  the corresponding unipotent root subgroup with respect to $T$. We denote by $E(\Phi,R)=E_\P(\Phi,R)$the elementary subgroup generated by all $X_\alpha$. 

Let $\Delta$ be a subsystem of $\Phi$. We denote by $E(\Delta,R)$ the subgroup of $G(\Phi,R)$  generated by all $X_\alpha$, where $\alpha\in \Delta$. It is called an (elementary) {\it subsystem subgroup}. It can be shown
that it is an elementary subgroup of a Chevalley group $G(\Delta,R)$ embedded into the group $G(\Phi,R)$. Here the lattice between $\euQ(\Delta)$ and $\P(\Delta)$ is an orthogonal projection of $\P$ to the corresponding subspace.

We are going to describe intermediate subgroups between $E(\Delta,R)$ and $G(\Phi,R)$. For
each case considered in the present paper, we fix a basic representation $V$ of the group $G(\Phi,R)$ (i.e., a representation that has a basis consisting of weight vectors such that
the Weyl group acts transitively on the set of nonzero weights, see \cite{atlas}) and we assume that the lattice $\P$ is such that this representation is faithful (therefore, $\P$ is uniquely
determined by the representation).

Here is the list of the cases we consider:

\renewcommand{\theenumi}{\alph{enumi}}
\begin{enumerate}

\item $\Phi=D_l$, $\Delta=A_{l-1}$ ($l\ge 5$), $V=V_{\varpi_l}$ (half-spin representation).

\item $\Phi=E_6$, $\Delta=D_5$, $V=V_{\varpi_l}$.

\item $\Phi=E_7$, $\Delta=E_6$, $V=V_{\varpi_7}$.
\end{enumerate}
\renewcommand{\theenumi}{\arabic{enumi}}

\begin{figure}[h]
\begin{center}

\includegraphics[scale=0.2]{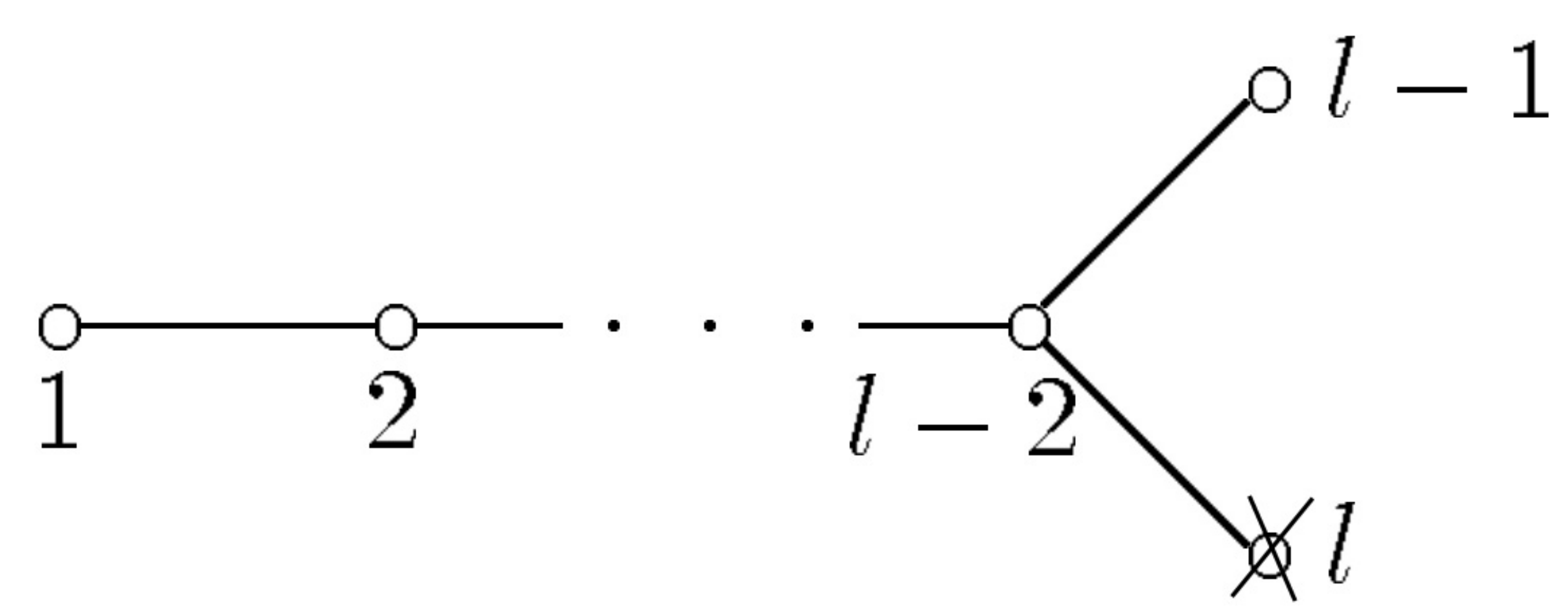}

\phan

\includegraphics[scale=0.2]{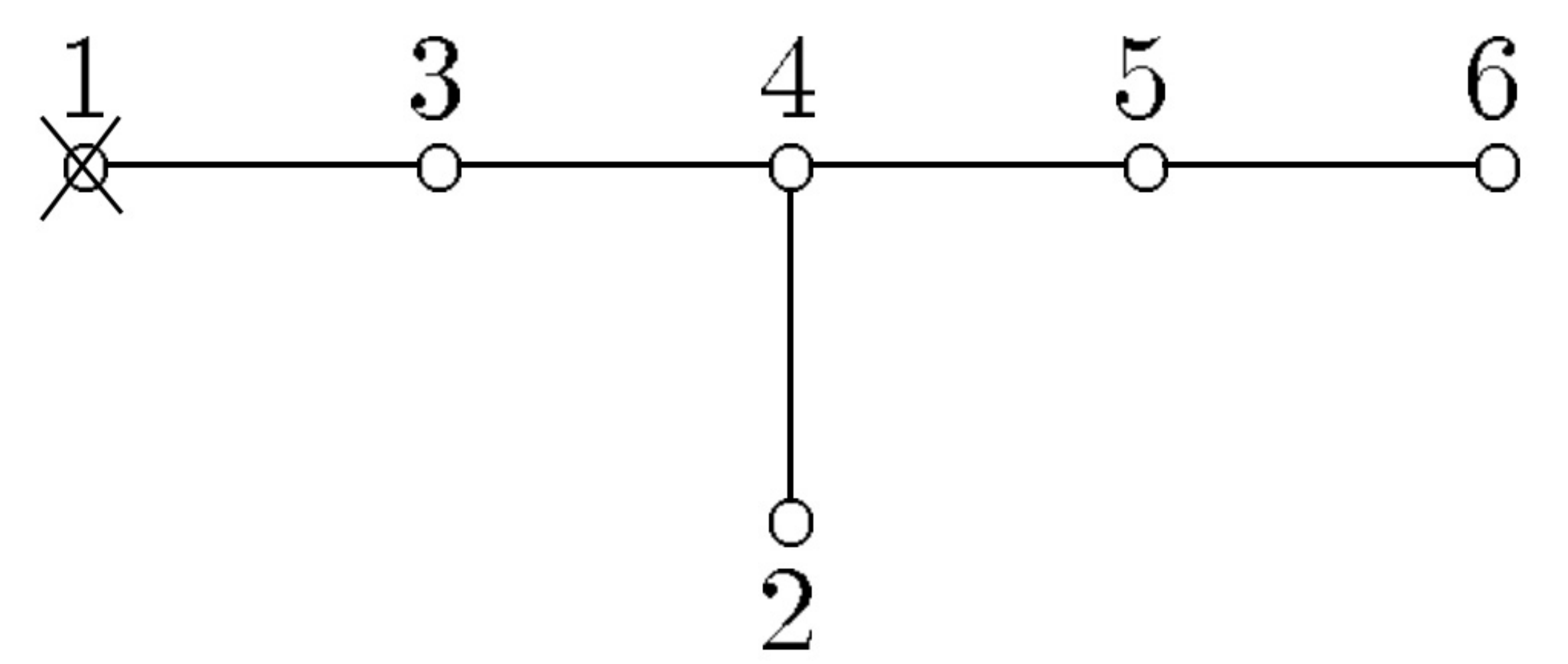}

\phan

\includegraphics[scale=0.2]{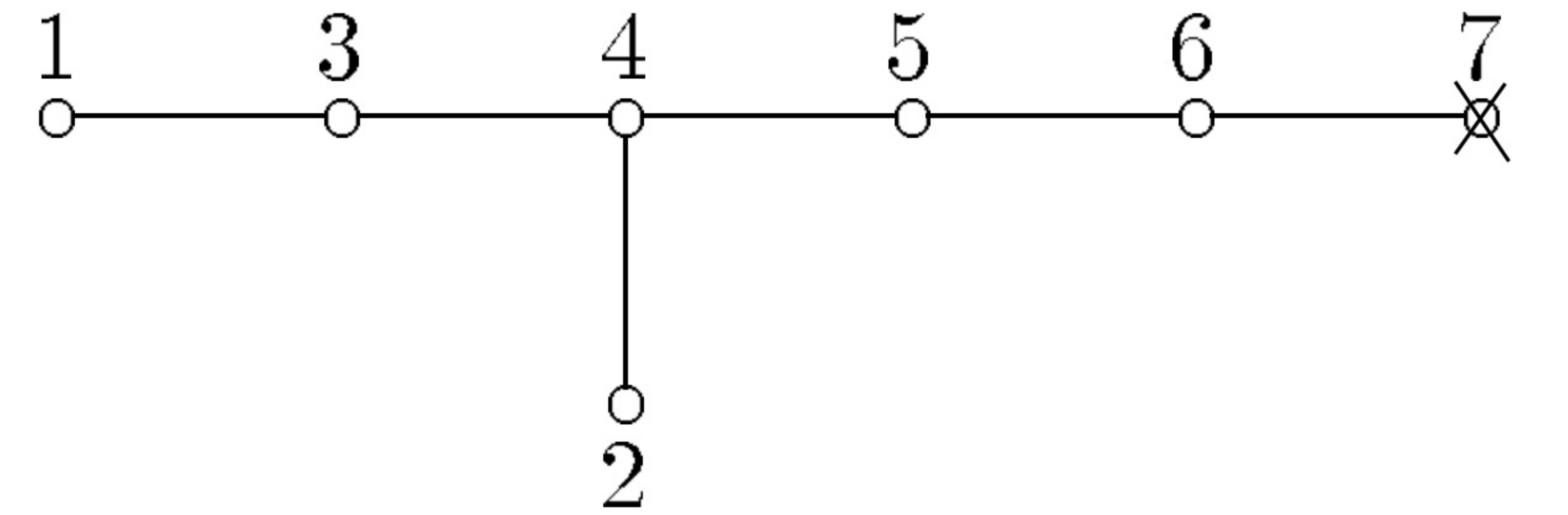}
\end{center}
\caption{}
\label{Dynkin}
\end{figure}

The embedding of $\Delta$ into $\Phi$ is obtained by crossing out a vertex of the Dynkin diagram
as shown in Figure~\ref{Dynkin}. 

Let us fix an order on the system $\Phi$  and denote by $\Pi=\{\alpha_1$, $\ldots$, $\alpha_l\}$ the corresponding
set of fundamental roots; its enumeration is also shown in Figure~\ref{Dynkin}. Let $\alpha^{(1)}$ be the root that corresponds to the vertex crossed out, and $\alpha^{(2)}$ the root that corresponds to the vertex that is adjacent to the crossed out one.

Let $\Delta'$ be the subsystem of $\Delta$ obtained by crossing out both $\alpha^{(1)}$ and $\alpha^{(2)}$. Furthermore, let $\Delta''$ be equal to $\Delta'$ in cases (b) and (c) and to its irreducible component distinct from $A_1$ (it is important that $l\ge 5)$ in case (a).

The symbols $\Omega^\pm$, $\Sigma_{\lambda_1}^{\pm,0}$, $\Delta_{\lambda_1}$, and $(\Delta\cap\Delta_{\lambda_1})'$ will be introduced in Lemmas \ref{Weylorbits} and \ref{SigmaStructure}.

\subsection{Affine schemes}

The functor $G(\Phi,-)$ from the category of rings to the category
of groups is an affine group scheme (a Chevalley--Demazure scheme). This means that
its composition with the forgetful functor to the category of sets is representable, i.e.,
$$
G(\Phi,R)=\Hom (\Z[G],R)\tp
$$
The ring $\Z[G]$ here is called the {\it ring of regular functions} on the scheme $G(\Phi,-)$. 

An element $g_{\gen}\in G(\Phi,\Z[G])$, that corresponds to the identity ring homomorphism is
called the {\it generic element} of the scheme  $G(\Phi,-)$. This element has a universal property:
for any ring $R$ and for any $g\in G(\Phi,R)$, there exists a unique ring homomorphism 
$$
\map{f}{\Z[G]}{R}
$$
such that $f_*(g_{\gen})=g$. For details about application of the method of generic element
to the problems similar to our problem see the paper of Alexei Stepanov \cite{StepUniloc}.  We will
use this method in the second part of Lemma \ref{A2} to avoid problems with zero divisors of
the ring $R$.

Later, we will define the scheme $\Root(-)$ of the root type elements and introduce a
similar notation for it.

The identity element of the group $G(\Phi,\Z)$ corresponds to the augmentation homomorphism
$$
\map{\eps}{\Z[G]}{\Z}\tp
$$
Its kernel is called the {\it augmentation} ideal. We denote it by $I_{\aug}$.

\subsection{Representation $V$}

Let $\Lambda$ be the set of weights of the representation $V$.  A basis
of $V$ V is indexed by the weights from $\Lambda$. We denote this basis by $\{v^\lambda\colon \lambda\in\Lambda\}$. 

The order we fixed on the system $\Phi$ induces a partial order on the set $\Lambda$:: a weight $\lambda$ is bigger than a weight  $\mu$ if $\lambda-\mu$ is the sum of positive roots.

Root elements of the group $G(\Phi,R)$ act on the elements of this basis as follows:
$$
x_\alpha(\xi)v^\lambda=
\begin{cases}
v^\lambda+v^{\lambda+\alpha} c_{\lambda\alpha}\xi \quad&\text{if } \lambda+\alpha\in\Lambda\tc\\
v^\lambda \quad&\text{if } \lambda+\alpha\notin\Lambda\tc
\end{cases}
$$
where the structure constants $c_{\lambda\alpha}$ of the action are equal to $\pm 1$ (Matsumoto’s lemma \cite{Matsumoto69}, \cite{Stein}). Moreover, we choose $\{v^\lambda\}$ to be a crystal basis, i.e., $c_{\lambda\alpha}=1$ if $\alpha\in\pm\Pi$. The
existence of such a basis was proved in \cite{VavThirdlook}, \cite{Borel}, \cite{VavSigns}.

Therefore, the weight diagram describes the action of $E(\Phi,R)$ on the module $V$ completely. This action has a unique extension that is an algebraic action of the group $G(\Phi,R)$.

Weights from $\Lambda$ correspond to the vertices of the weight diagram (see \cite{atlas}). To describe
the action of the group $G(\Delta,R)$ on the module $V$, one should remove the edges of the
weight diagram labeled by $\alpha^{(1)}$. The new diagram will have a number of connected
components, which can be enumerated $\Lambda_0$, $\Lambda_1$, $\ldots$ in such a way that whenever $i>j$ any weight from $\Lambda_i$ is less that any weight from $\Lambda_j$ with respect to the order we fix. The
set $\Lambda_0$ contains only one weight --- the highest one, we denote it by $\lambda_0$.

For any pair of weights $\lambda$,$\mu\in\Lambda$ we denote by $d(\lambda,\mu)$ the distance between them (i.e., the length of the shortest path) in the weight graph (two weights are adjacent in the weight graph if their difference belongs to $\Phi$).

\subsection{Group theoretic notation}
\begin{itemize}
\item Recall that for abstract groups $A,B\le G$, the transporter from $A$ to $B$  is the set
$$
\Tran_G(A,B)=\{g\in G\colon gAg^{-1}\sub B\}\tp
$$

\item We denote by $N_G(\Gamma)$  the normaliser of the group $\Gamma$  in the group $G$. 

\item If the group $G$ acts on the set $X$ and $x\in X$, we denote by $\Stab_G(x)$ the stabiliser of the element $x$.

\item Commutators are left normalised:
$$
[x,y]=xyx^{-1}y^{-1}\tp
$$

\item If $X$ is a subset of the group $G$, we denote by $\<X\>$ the subgroup generated by $X$.

\item We denote by $D^i G$ the $i$th member of the derived series, i.e., $D^0 G=G$ and $D^{i+1}G=[D^iG,D^iG]$.
\end{itemize}

\subsection{Matrices}
We identify any element $g$ of the group  $G(\Phi,R)$ with the matrix of its
action on the module $V$ in the crystal basis. Rows and columns of such a matrix are
indexed by weights from $\Lambda$. We denote by $g_{\lambda,\mu}$ the entry of this matrix in the $\lambda$th row and the $\mu$th column. We denote by $g_{*,\mu}$ the $\mu$th column of the matrix $g$, which we identify with a vector from $V$. Similarly, we denote by $g_{\lambda,*}$ the $\lambda$th row of the matrix $g$,  which we identify with a covector from  $V^*$.

\subsection{Relative subgroups}

For a root system $\Phi$, a ring $R$ and an ideal $I\unlhd R$  we denote by $\rho_I$ the reduction homomorphism
$$
\map{\rho_I}{G(\Phi,R)}{G(\Phi,\fact{R}{I})}\tc
$$
induced by the projection of $R$ to the quotient ring $\fact{R}{I}$.
We denote by $G(\Phi,R,I)$ the principal congruence subgroup of the group $G(\Phi,R)$
$$
G(\Phi,R,I)=\Ker \rho_I\tp
$$
We denote by $CG(\Phi,R,I)$ the full congruence subgroup, i.e., the inverse image of the
center of the group $G(\Phi,\fact{R}{I})$ under the reduction homomorphism $\rho_I$.

We denote by $E(\Phi,I)$ the subgroup generated by the root elements of the level $I$:
$$
E(\Phi,I)=\<\{x_\alpha(\xi)\colon \alpha\in \Phi, \gap \xi\in I\}\>\tp
$$

We denote by $E(\Phi,R,I)$ the relative elementary subgroup of the level $I$, i.e., the
normal closure of the group $E(\Phi,I)$ in the group $E(\Phi, R)$. As a subgroup, $E(\Phi,R,I)$ is generated be the elements $z_\alpha(\xi,\zeta)=x_\alpha(\zeta)x_{-\alpha}(\xi)x_\alpha(-\zeta)$, where $\alpha\in \Phi$, $\xi\in I$ and $\zeta\in R$ (see \cite{TitsCongruence}).

The following facts about this group are well known.

{\lem\label{standcomm}{\rm (}Standard commutator formula {\rm )} If $\rk\Phi\ge 2$, then $E(\Phi,R,I)$ is a normal subgroup of $G(\Phi,R)$ and
$$
[E(\Phi,R),CG(\Phi,R,I)]\le E(\Phi,R,I)\tp
$$}

{\lem\label{ISquare}{\rm (}see \cite{TitsCongruence}{\rm ) } Let $I\unlhd R$. Then
$$
E(\Phi,R,I^2)\le E(\Phi,I)\tp
$$}

\subsection{Nilpotent structure of $K_1$}

We will use the following fact proved in \cite{BHVstrike}.

{\lem\label{K1nilp} For an arbitrary commutative ring $R$ of finite Bass--Serre dimension {\rm (}in particular, for any finitely generated ring{\rm )} and an arbitrary idea $I\unlhd R$, the quotient group $\fact{G(\Phi,R,I)}{E(\Phi,R,I)}$ has a nilpotent normal subgroup with Abelian quotient group {\rm (}in particular, it is solvable{\rm )}.}

\subsection{Parabolic subgroups}

Assume that we consider one of cases (a)–(c).

We denote by $P=P_{\alpha^{(1)}}$ the parabolic subgroup of the group $G(\Phi,R)$,  that corresponds to the root $\alpha^{(1)}$.

This subgroup $P$ coincides with the stabiliser of the line generated by the highest weight vector. The stabiliser of the line generated by the vector of the weight $\lambda\in \Lambda$ is also a parabolic subgroup (it is conjugate to $P$ by an element of the Weyl group\footnote{We are considering basic representation without zero weight. This means exactly that the Weyl group acts transitively on the set of weights.}, that maps the highest weight to the weight $\lambda$), we will denote it by $P_\lambda$.

We denote the corresponding opposite parabolic subgroups by $P^-$ and $P_\lambda^-$, respectively. The subgroup $P_\lambda^-$ can be described as the stabiliser of the line generated by the covector that corresponds to the weigh $\lambda$.
In other words, the subgroup $P_\lambda$ consist of the matrices such that their $\lambda$th  column is a multiple of the corresponding column of the identity matrix. Similarly, the subgroup $P_\lambda^-$ consists of the matrices such that their $\lambda$th row is a multiple the corresponding row of the identity matrix.
 
The unipotent radicals of these parabolic subgroups will be denoted by $U$, $U_\lambda$, $U^-$, and $U^-_\lambda$, respectively. The correspondent Levi subgroups will be denoted by $L$ and $L_\lambda$ (opposite parabolic subgroups have the same Levi subgroups).

The following statements are well known.

{\lem \label{parabolic} In each of cases {\rm (}a{\rm )}--{\rm (}c{\rm )}
\begin{enumerate}
\item
$$
L_\lambda=P_\lambda\cap P^-_\lambda\tp
$$

\item 
$$
P_\lambda=U_\lambda\leftthreetimes L_\lambda\quad\quad\text{and}\quad\quad P^-_\lambda=U_\lambda^-\leftthreetimes L_\lambda\tp
$$

\item
$$
U_\lambda=\<\{x_\alpha(\xi)\colon \xi\in R,\gap \lambda-\alpha\in\Lambda\}\>\tc
$$
$$
U_\lambda^-=\<\{x_\alpha(\xi)\colon \xi\in R,\gap \lambda+\alpha\in\Lambda\}\>\tp
$$

\item The groups $U_\lambda$ and $U^-_\lambda$ are Abelian.
\end{enumerate}} 
 
\section{Combinatorial lemmas} 

In this section, we assume that we are considering one of cases (a)--(c), which were
listed at the beginning of the paper.

Let $W(\Phi)$ and $W(\Delta)$ be the Weyl groups of the corresponding root systems.
 
{\lem\label{alphaonecoefficient} The coefficient of $\alpha^{(1)}$  in the decomposition of the maximal root is equal to one. The coefficient of $\alpha^{(2)}$ is equal to two.}
\begin{proof}
Direct calculation in each case.
\end{proof} 

{\lem\label{Weylorbits} The set $\Phi\sm\Delta$ has exactly two $W(\Delta)$-orbits$:$ the roots such that in their expansion in fundamental roots the coefficient of $\alpha^{(1)}$ is equal to $1$, and the roots such that it is equal to $-1$. We will denote these orbits by $\Omega^+$ and $\Omega^-$ respectively.}
\begin{proof}
We prove that if a root $\alpha\in\Phi$, has coefficient of $\alpha^{(1)}$ equal to $1$,then it lies in one $W(\Delta)$--orbit with the maximal root (the proof for the opposite roots is the same).

One can obtain the maximal root from the root $\alpha$ by adding fundamental roots step
by step. All these fundamental roots are distinct from $\alpha^{(1)}$, (hence they belong to $\Delta$) because the coefficient of $\alpha^{(1)}$ in the expansion of $\alpha$ is already maximal possible. If one root is obtained from another by adding a root from $\Delta$, then it is also obtained by reflection with respect to this root (because all roots have the same length). Therefore, all roots in the chain from  $\alpha$ to the maximal root lie in the same $W(\Delta)$-orbit.
\end{proof}
 
{\lem \label{WeylorbitsinDelta}
  The set $\Delta\sm\Delta'$ has exactly two $W(\Delta')$-orbits$:$ the roots such that in their
  expansion in fundamental roots the coefficient of $\alpha^{(2)}$ is equal to $1$, and the roots such that it is equal to $-1$.}
  \begin{proof}
The arguments are similar to the proof of the previous lemma.
\end{proof}   
 
 {\lem \label{RootInDeltaComb} Suppose $\beta\in\Phi\sm\Delta$; then there exists a root $\alpha\in \Delta$ such that $\alpha+\beta\in \Phi$.}
 \begin{proof}
 Using the action of the Weyl group $W(\Delta)$ we may assume that $\beta=\pm\alpha^{(1)}$. In
 this case, we can take $\alpha=\pm\alpha^{(2)}$.
 \end{proof}
 
 {\lem\label{OmegaAndHighestWeight}
Suppose $\alpha\in \Omega^+$; then
$
\lambda_0-\alpha\in\Lambda\tp
$  
  }
  \begin{proof}
  Note that the weight $\lambda_0$ is fixed by the group $W(\Delta)$. Hence using the action of
  the Weyl group $W(\Delta)$ we may assume that $\alpha=\alpha^{(1)}$. In each of cases (a)--(c), one can see on the weight diagram that the highest weight is, indeed, adjacent to another one by the edge labeled with $\alpha^{(1)}$. 
  \end{proof}

{\lem\label{notroot} Suppose $\alpha\in\Phi$ and $\lambda$, $\rho$, $\lambda-\alpha$, $\rho+\alpha\in\Lambda$ are such that $\lambda-\alpha\ne \rho$; then $d(\lambda,\rho)\ge 2$.}
\begin{proof}
First, using the action of the group $W(\Phi)$, we may assume that $\lambda=\lambda_0$ is the
highest weight. Second, using the action of the group $W(\Delta)$, we may also assume that $\alpha=\alpha^{(1)}$. 

Consider the difference $\lambda-(\rho+\alpha)$. To calculate this difference, one should consider the path on the weight diagram between the vertices $\lambda$ and $\rho+\alpha$ (since $\lambda$ is the highest weight, this path can always be chosen going from the left to the right), and take the sum of the fundamental roots corresponding to the edges of this path. Note that $\alpha^{(1)}$ appears in this sum at least once because this is the only edge incident to the vertex $\lambda_0$, and $\lambda\ne \rho+\alpha$ by assumption. Therefore, in the expansion of $\lambda-\rho$ in fundamental roots, the coefficient of$\alpha^{(1)}$  is at least 2; hence by Lemma~\ref{alphaonecoefficient}, $\lambda-\rho$ is not a root.
\end{proof}  

{\lem\label{SigmaStructure}
Suppose $\lambda_1\in\Lambda_1$. We introduce the following notation: 
$
\Sigma_{\lambda_1}=\{\alpha\in\Phi\colon \lambda_1-\alpha\in\Lambda\}\tp
$
Next, we split this set in three disjoint subsets:
$$
\Sigma_{\lambda_1}=\Sigma_{\lambda_1}^-\sqcup \Sigma_{\lambda_1}^0\sqcup \Sigma_{\lambda_1}^+\tc
$$
where
$$
\Sigma_{\lambda_1}^-=\Sigma_{\lambda_1}\cap\Omega^-\tc\quad\Sigma_{\lambda_1}^0=\Sigma_{\lambda_1}\cap\Delta\tc\quad \Sigma_{\lambda_1}^+=\Sigma_{\lambda_1}\cap\Omega^+\tp
$$

Let $\Delta_{\lambda_1}$ be the image of $\Delta$ under reflection with respect to $\lambda_0-\lambda_1$. Then the following holds.
\begin{enumerate}
\item $\Sigma_{\lambda_1}^-=\{\lambda_1-\lambda_0\}$.

\item $\Sigma_{\lambda_1}^0\ne\emp$ and for any $\alpha\in\Sigma_{\lambda_1}^0$ there exists a root $\beta\in\Sigma_{\lambda_1}^0$ such that $\alpha-\beta\in\Phi$.

\item  For any $\beta\in \Sigma_{\lambda_1}^+$ there exists a root $\gamma\in \Delta\cap\Delta_{\lambda_1}$ such that $\beta+\gamma\in\Phi$, but \linebreak $\lambda_1-\lambda_0+\gamma\notin\Phi$. 

\item The subsystem $\Delta\cap\Delta_{\lambda_1}$ has an irreducible component $(\Delta\cap\Delta_{\lambda_1})'$ distinct from $A_1$. Next, for any weight $\mu\in \Lambda_1\sm\{\lambda_1\}$ there exists a weight $\nu\in \Lambda_1\sm\{\lambda_1\}$ such
that $\mu-\nu\in(\Delta\cap \Delta_{\lambda_1})'$.

\item For any $\alpha\in\Omega^+$ such that $\<\alpha,\lambda_0-\lambda_1\>=1$, there exists a root $\gamma\in\Sigma_{\lambda_1}^0$ such that $\alpha+\gamma\in\Phi$.
\end{enumerate}
}
\begin{proof}
Using the action of the group $W(\Delta)$, we may assume that $\lambda_0-\lambda_1=\alpha^{(1)}$. Therefore, we have $\Delta\cap\Delta_{\lambda_1}=\Delta'$, and in Item (4) we have $(\Delta\cap\Delta_{\lambda_1})'=\Delta''$. Let $w_{\alpha^{(1)}}\in W(\Phi)$ be reflection with respect to the root $\alpha^{(1)}$. This reflection maps $\lambda_0$ to $\lambda_1$; hence it maps $\Omega^+$ to $\Sigma_{\lambda_1}$. 

We prove items (1) and (2). Since $w_{\alpha^{(1)}}\alpha^{(1)}=-\alpha^{(1)}$, we have $-\alpha^{(1)}\in \Sigma_{\lambda_1}^-$.
Suppose that $\alpha\in\Omega^+$ and $\alpha\ne\alpha^{(1)}$, then $\<\alpha,\alpha^{(1)}\>$ is equal either to 0 or to 1 (because $\alpha+\alpha^{(1)}\notin\Phi$). In the first case, we have
$$
w_{\alpha^{(1)}}\alpha=\alpha\in\Sigma_{\lambda_1}^+\tp
$$ 
In the second case, we have
$$
w_{\alpha^{(1)}}\alpha=\alpha-\alpha^{(1)}\in\Sigma_{\lambda_1}^0\tp
$$
Item (1) is, therefore, proved.

Furthermore, if $w_{\alpha^{(1)}}\alpha\in \Sigma^0_{\lambda_1}$, then $\<w_{\alpha^{(1)}}\alpha,\alpha^{(1)}\>=-1$ (because $\<\alpha,\alpha^{(1)}\>=1$).Hence
the coefficient of $\alpha^{(2)}$ in the expansion of $w_{\alpha^{(1)}}\alpha$ is equal to 1. Conversely, if the coefficient of $\alpha^{(2)}$ in the expansion of the root $w_{\alpha^{(1)}}\alpha\in\Delta$ is equal to 1, then
$$
\alpha=w_{\alpha^{(1)}}w_{\alpha^{(1)}}\alpha=w_{\alpha^{(1)}}\alpha+\alpha^{(1)}\in\Omega^+\tp
$$
Hence $w_{\alpha^{(1)}}\alpha\in\Sigma_{\lambda_1}^0$.

Therefore, the set $\Sigma_{\lambda_1}^0\sub\Delta$ is an analog of the set $\Omega^+$  for the pair $(\Delta,\Delta')$. Now,
clearly, it is not empty, and the second part of item (2) follows from an analog of Lemma~\ref{RootInDeltaComb} for the pair $(\Delta,\Delta')$.

We prove item (3). From the above, it follows that $\<\beta,\alpha^{(1)}\>=0$. Hence the coefficient
of $\alpha^{(2)}$ in the expansion of $\beta$ ? is equal to 2. Therefore, the maximal root can be obtained from $\beta$ by adding (= reflecting with respect to) simple roots from $\Delta'$ step by step. Using the action of the group $W(\Delta')$ we may assume that $\beta$ is a maximal root. Now we can take $\gamma=-\alpha_k$, where $\alpha_k$ is a simple root that corresponds to the vertex of the Dynkin diagram that is adjacent to the additional vertex in the affine diagram.\footnote{It is important that $l\ge 5$, and that we do not consider the case $D_{l-1}\le D_l$. Theorem \ref{sandwich} is not true for these subsystems because we can take $H=G(B_{l-1},R)$.}

The first part of Item (4) has already been proved above. Moreover, $(\Delta\cap\Delta_{\lambda_1})'=\Delta''$. Looking at the weight diagram, one can easily check the second part.

Let us prove Item (5). By hypothesis, $\alpha-\alpha^{(1)}\in\Delta$, and coefficient of $\alpha^{(2)}$ in its decomposition is equal to one. Using the action of the group $W(\Delta')$, we can assume that the sum of coefficients in the decomposition of $\alpha$ in simple roots is maximal among all such sums for the roots from the $W(\Delta')$-orbit of $\alpha$. There is a simple root $\alpha_k$, such that $\alpha+\alpha_k\in\Phi$. The root $\alpha_k$ cannot be equal to $\alpha^{(1)}$ because the coefficient of $\alpha^{(1)}$ is
already equal to 1. Also the root $\alpha_k$ cannot be in $\Delta'$;  otherwise, $\alpha+\alpha_k$ would be a root from the $W(\Delta')$-orbit of $\alpha$ with bigger sum of coefficients. Hence $\alpha+\alpha^{(2)}\in\Phi$ and we can take $\gamma=\alpha^{(2)}$.
\end{proof}

{\lem\label{weights}
Let $\lambda_1\in\Lambda_1$; then:
\begin{enumerate}
\item there exists $\mu\in\Lambda_1$ such that $d(\lambda_1,\mu)=1$;

\item let $\nu\in\Lambda_1$ be a weight such that $d(\lambda_1,\nu)=1$;  then there exists $\mu\in\Lambda_1$ such that $d(\mu,\nu)= d(\lambda_1,\mu)=1$.
\end{enumerate}
}
\begin{proof}
To each weight $\mu\in \Lambda_1$ such that $d(\lambda_1,\mu)=1$,  we assign a root $\lambda_1-\mu\in \Sigma_{\lambda_1}^0$. This correspondence is bijective. Therefore, this lemma is merely a restatement of Item~2 of Lemma~\ref{SigmaStructure}.
\end{proof}

\section{Level computation} 
 
 For now, let $\Delta$, $\Phi$, and $\P$ be arbitrary. Let $H$ be an overgroup of the group $E(\Delta,R)$, i.e.
 $
 E(\Delta,R)\le H\le G(\Phi,R)\tp
 $
 For each root $\alpha\in\Phi\sm \Delta$ we set
 $
 I_\alpha=\{\xi\in R\colon x_\alpha(\xi)\in H\}\tp
 $
 In the paper \cite{VSch} (Lemmas 1 and 3), the following was proved.
 {\lem\label{level}
 \begin{enumerate}
 \item The set $I_\alpha$ is an ideal of the ring $R$.
 
 \item The ideal $I_\alpha$ depends only on the $W(\Delta)$-orbit of the root $\alpha$.
 \end{enumerate}}
 
Now we assume that we are considering one of cases (a)--(c) that were listed in the
begining of the paper. Recall that by Lemma~\ref{Weylorbits}, the set $\Phi\sm\Delta$ has exactly two $W(\Delta)$-orbits: $\Omega^+$ and $\Omega^-$. We denote the corresponding ideals $I_\alpha$ by $I^+$ and $I^-$ respectively.

The pair $\sigma=(I^+,I^-)$ is called the {\it level} of the overgroup $H$. We will write $\sigma=\lev(H)$.

Now let $\sigma=(I^+,I^-)$ be an arbitrary pair of ideals of the ring $R$, We introduce the
following notation:
\begin{gather*}
E(\Phi,\Delta,R,\sigma)=\<x_\alpha(\xi)\colon \alpha\in\Delta,\,\xi\in R\;\\ \text{or}\; \alpha\in\Omega^+,\, \xi\in I^+\;\text{or}\; \alpha\in\Omega^-,\, \xi\in I^-\>\le G(\Phi,R)\tc\\
G(\Phi,\Delta,R,(R,I^-))=\left(\rho_{I^-}\right)^{-1}(P)\le G(\Phi,R)\tc\\
G(\Phi,\Delta,R,(I^+,R))=\left(\rho_{I^+}\right)^{-1}(P^-)\le G(\Phi,R)\tc\\
G(\Phi,\Delta,R,\sigma)=G(\Phi,\Delta,R,(R,I^-))\cap G(\Phi,\Delta,R,(I^+,R))\le G(\Phi,R)\tp
\end{gather*}
It is easily seen that
 $$
 E(\Phi,\Delta,R,\sigma)\le G(\Phi,\Delta,R,\sigma)\tc
 $$
and also that
 $$
 \sigma\le\lev\left( E(\Phi,\Delta,R,\sigma)\right)\le\lev\left(G(\Phi,\Delta,R,\sigma)\right)\le\sigma
 $$
 (we write $\sigma_1\le\sigma_2$ if $I^+_1\sub I^+_2$ and $I^-_1\sub I^-_2$). Hence in fact, we have the following lemma.
 {\lem\label{LevelsOfEandG} We have
 $$
 \lev\left(E(\Phi,\Delta,R,\sigma)\right)=\lev\left(G(\Phi,\Delta,R,\sigma)\right)=\sigma\tp
 $$}

Therefore, any pair of ideals is a level of some overgroup of the group $E(\Delta,R)$. Let
us also compute the level of the group $E(\Phi,\Delta,R,\sigma)$.
{\lem \label{NormalizerLevel}. We have
$$
\lev(N_{G(\Phi,R)}(E(\Phi,\Delta,R,\sigma)))=\sigma\tp
$$}
\begin{proof}
First,
$$
\lev(N_{G(\Phi,R)}(E(\Phi,\Delta,R,\sigma)))\ge\lev(E(\Phi,\Delta,R,\sigma))=\sigma\tp
$$

We prove the inverse inclusion. Assume that
$
x_\beta(\xi)\in N_{G(\Phi,R)}(E(\Phi,\Delta,R,\sigma))\tc
$

where $\beta\in\Omega^+$, but $\xi\notin I^+$ (the proof for $\Omega^-$ is similar). Take $\alpha$ from Lemma~\ref{RootInDeltaComb}; then
$
x_\alpha(1)\in E(\Delta,R)\le E(\Phi,\Delta,R,\sigma)\tp
$
Hence we have
$
[x_\alpha(1),x_\beta(\xi)]\in E(\Phi,\Delta,R,\sigma)\tp
$
On the other hand, we have
$
[x_\alpha(1),x_\beta(\xi)]=x_{\alpha+\beta}(\pm\xi)\tc
$
where $\alpha+\beta\in\Omega^+$, which contradicts the fact that 
$
\lev(E(\Phi,\Delta,R,\sigma))=\sigma\tp
$
\end{proof}

\section{The statement of the main result}

In the present paper, we prove the following theorem

{\thm\label{sandwich} Assume that we are considering one of cases (a)–(c) listed at the begining of the paper, i.e
\renewcommand{\theenumi}{\alph{enumi}}
\begin{enumerate}

\item $\Phi=D_l$, $\Delta=A_{l-1}$ (where $l\ge 5$), $V=V_{\varpi_l}$ {\rm (}half-spin representation{\rm )}.

\item $\Phi=E_6$, $\Delta=D_5$, $V=V_{\varpi_l}$.

\item $\Phi=E_7$, $\Delta=E_6$, $V=V_{\varpi_7}$.
\end{enumerate}

Let $R$ be a commutative ring {\rm(}associative with unit{\rm )}. . Then for any group $H$  such that
$$
E(\Delta,R)\le H \le G(\Phi,R)\tc
$$
there exists a unique pair of ideals $\sigma=(I^+,I^-)$ of the ring $R$, such that
$$
E(\Phi,\Delta,R,\sigma)\le H\le N_{G(\Phi,R)}\left(E(\Phi,\Delta,R,\sigma)\right)\tp
$$}

Note that uniqueness has already been proved: indeed, by Lemmas~\ref{LevelsOfEandG} and~ \ref{NormalizerLevel} the pair $\sigma$ must be equal to the pair $\lev(H)$. If we set $\sigma=\lev(H)$, then the left inclusion in the theorem holds true by the definition of a level. Thus, it remains to prove the right
inclusion.
 
It is natural to call the subgroup $E(\Phi,\Delta,R,(R,0))$ the elementary parabolic subgroup.
Note that Theorem~\ref{sandwich} gives us an overgroup description for such a subgroup.

{\cor  Assume that we are considering one of cases (a)–(c) that are listed at the begining of the paper. Let $R$ be a commutative ring. Then for any group $H$ such that
$$
E(\Phi,\Delta,R,(R,0))\le H \le G(\Phi,R)\tc
$$
there exists a unique ideal $I^-$ of $R$ such that
$$
E(\Phi,\Delta,R,(R,I^-))\le H\le N_{G(\Phi,R)}\left(E(\Phi,\Delta,R,(R,I^-))\right)\tp
$$
}

Similarly, we obtain the same result for the subgroup $E(\Phi,\Delta,R,(0,R))$.

To formulate the second main result, we divide cases (a)–(c) into two groups. Let $n+1$ --- be the number of connected components $\Lambda_0$,$\ldots$,$\Lambda_n$ of the weight diagram after removing the edges labeled with $\alpha^{(1)}$. Note that for $i\ne 0$,$n$ the component $\Lambda_i$ has at least two weights. Indeed, given a weight from such a component we can always add a simple root to it, and we can always subtract a simple root from it (so that the result is still in $\Lambda$). At least one of these roots is distinct from $\alpha^{(1)}$ because the representation is minuscule.

Therefore, two cases are possible. We say that we are considering the case of the {\it first
type} if $\Lambda_0$ is the only component containing only one weight (that includes case (a) with $l$ odd, and also case (b)). We say that we are considering the case of the second type if
there are two such components, $\Lambda_0$ and $\Lambda_n$ (that includes case (a) with $l$ even, and also case (c)). In the second case it is easily seen that $\Lambda_n=\{-\lambda_0\}$.

{\thm\label{Normalizer} Let $\sigma=(I^+,I^-)$ be a pair of ideals of the ring $R$; then:
 \begin{enumerate}
\item for the cases of the first type, we have
$$
N_{G(\Phi,R)}\left(E(\Phi,\Delta,R,\sigma)\right)=G(\Phi,\Delta,R,\sigma)\tp
$$

\item  for the cases of the second type, the group $N_{G(\Phi,R)}\left(E(\Phi,\Delta,R,\sigma)\right)$ consists of
exactly the elements $g$ of the group $G(\Phi,R)$ that satisfy the following conditions:
\begin{align*}
& g_{\lambda_0,\lambda}\in I^+\quad\forall\lambda\in \Lambda\sm\{\lambda_0,-\lambda_{0}\}\tc\\
& g_{\lambda_0,-\lambda_0}I^-\sub I^+\tc\\
& (g^{-1})_{\lambda,\lambda_0}\in I^-\quad\forall\lambda\in \Lambda\sm\{\lambda_0,-\lambda_{0}\}\tc\\
& (g^{-1})_{-\lambda_0,\lambda_0}I^+\sub I^-\tp
\end{align*} 
\end{enumerate}}

\section*{Agreement}

At this point, we begin the proof of Theorem~\ref{sandwich}, i.e., below we always assume that we are considering one of cases (a)--(c) that were listed at the begining of the paper

\section{Normalizer of the subgroup $E(\Phi,\Delta,R,\sigma)$} 
 
 {\lem\label{UcapGRI} Let $I\unlhd R$; then
 $$
 U\cap G(\Phi,R,I)=\<\{x_\alpha(\xi)\colon \xi\in I,\gap \alpha\in \Omega^+\}\>\tc
 $$
and the same is true for $U^-$.}
 \begin{proof}
 Obviously, the right-hand side is contained in the left-hand side, let us prove the
 inverse inclusion. Let $g\in  U\cap G(\Phi,R,I)$. Since $g\in U$, by Items 3 and 4 of Lemma~\ref{parabolic}, the element $g$ can be writen as the product 
 $$
 g=\prod_{\alpha\in\Omega^+}x_\alpha(\xi_\alpha)\tc\qquad\text{where } \xi_\alpha\in R\tp
 $$ 
 Fixing $\ovl{\alpha}\in \Omega^+$, we prove that $\xi_{\ovl{\alpha}}\in I$. Let $\rho=\lambda_0-\ovl{\alpha}$; then by Lemma~\ref{OmegaAndHighestWeight}, we have $\rho\in\Lambda$. Note that for an arbitrary matrix $h$ and for $\alpha\in \Omega^+\sm \{\ovl{\alpha}\}$ we have
 $$
 (x_\alpha(\xi_\alpha)h)_{\lambda_0,\rho}=h_{\lambda_0,\rho}\tp
 $$
 Indeed, otherwise we have $\rho+\alpha\in\Lambda$, and by Lemma~\ref{OmegaAndHighestWeight} we also have $\lambda_0-\alpha\in\Lambda$, which contradicts Lemma~\ref{notroot} .
 
 Therefore,
 $$
 \xi_{\ovl{\alpha}}=\pm(x_{\ovl{\alpha}}(\xi_{\ovl{\alpha}}))_{\lambda_0,\rho}=\pm g_{\lambda_0,\rho}\in I\tp
 $$ 
 \end{proof}
 
 In the lemmas below, let $\sigma=(I^+,I^-)$, $I^\pm\unlhd R$, and let $J=(I^+\cap I^-)^2$.
 
 {\lem\label{ChevMats} Assume that $J=0$; then
 \begin{enumerate}
\item
 $
 G(\Phi,\Delta,R,\sigma)=E(\Phi,\Delta,R,\sigma)T(\Phi,R)G(\Delta,R)\tp
 $
 \item the subgroup  $E(\Phi,\Delta,R,\sigma)$ is normal in $G(\Phi,\Delta,R,\sigma)$.
 \end{enumerate}} 
 \begin{proof}
First, note that the subgroup $T(\Phi,R)G(\Delta,R)$ normalises the subgroup $E(\Phi,\Delta,R,\sigma)$; hence Item 2 follows from Item 1. Indeed, $T(\Phi,R)$ normalises each of the subgroups $\{x_\alpha(\xi)\colon \xi\in I\}$, where $\alpha\in \Phi$ and $I$ equals $I^\pm$ or $R$. Furthermore, $G(\Delta,R)$ normalises $E(\Delta,R)$ (this is a special case of Lemma~\ref{standcomm}). Finally, if $\alpha\in \Omega^+$, $\xi\in I^+$ and $g\in G(\Delta,R)\le L$, then by Lemma~\ref{UcapGRI} we have
$$
[x_\alpha(\xi),g]\in U\cap G(\Phi,R,I^+)=\<\{x_\alpha(\xi)\colon \xi\in I^+,\gap \alpha\in \Omega^+\}\>\sub E(\Phi,\Delta,R,\sigma)\tp
$$
The case of $\Omega^-$ and $I^-$ is similar.

Therefore, it remains to prove Item 1. Obviously, the right-hand side is contained in the left-hand side, let us prove the inverse inclusion.
 
Let $g\in G(\Phi,\Delta,R,\sigma)$. Since all the rows and columns of an invertible matrix are unimodular, and all the entries in the first column except $g_{\lambda_0,\lambda_0}$ belong to $I^-$, the element $g_{\lambda_0,\lambda_0}$ is invertible modulo $I^-$. Similarly, we deduce that it is invertible modulo $I^+$. Hence it is invertible modulo  $J$, i.e., it is simply an invertible element of the ring $R$.
 
Therefore, we can apply to the matrix $g$ the special case of the Chevalley--Matsumoto
decomposition (see \cite{ChevalleySemiSimp}, \cite{Matsumoto69}, \cite{Stein}) asserting that if $g\in G(\Phi,R)$ and $g_{\lambda_0,\lambda_0}\in R^*$, then
 $$
 g=vg_1u\tc
 $$
 where $u\in U$, $v\in U^-$ and $g_1\in T(\Phi,R)G(\Delta,R)$.  Arguing as in the proof of Lemma~\ref{UcapGRI}, and using the fact that $g\in G(\Phi,\Delta,R,\sigma)$, we see that actually
 $$
 u\in\<\{x_\alpha(\xi)\colon \xi\in I^+,\gap \alpha\in \Omega^+\}\>\le E(\Phi,\Delta,R,\sigma)\tc
 $$
 and
  $$
 v\in\<\{x_\alpha(\xi)\colon \xi\in I^-,\gap \alpha\in \Omega^-\}\>\le E(\Phi,\Delta,R,\sigma)\tp
 $$
 To finish the proof, it remains to note that the subgroup $T(\Phi,R)G(\Delta,R)$ normalises the
 subgroup $E(\Phi,\Delta,R,\sigma)$.
 \end{proof}
 
 {\lem\label{DerivSer} If the ring $R$ is finitely generated, then there exists a natural number $N$ such that
 $$
 D^N G(\Phi,\Delta,R,\sigma)=E(\Phi,\Delta,R,\sigma)\tp
 $$}
 \begin{proof}
 First, the right-hand side is contained in the left-hand side for every $N$ because
 the group $E(\Phi,\Delta,R,\sigma)$ is obviously perfect. We prove the inverse inclusion.
 
We start with the case where $J=0$ and prove that
 $$
 D^i(G(\Phi,\Delta,R,\sigma))\le E(\Phi,\Delta,R,\sigma)D^i(T(\Phi,R)G(\Delta,R))\tp
 $$ 
 The proof is by induction on $i$. The base of induction, for $i=0$, 0, follows from Lemma~\ref{ChevMats}. Now, we pass to the induction step. By the inductive hypothesis, we have
 \begin{align*}
 &D^{i+1}(G(\Phi,\Delta,R,\sigma))\le\\ &\le [E(\Phi,\Delta,R,\sigma)D^i(T(\Phi,R)G(\Delta,R)),E(\Phi,\Delta,R,\sigma)D^i(T(\Phi,R)G(\Delta,R))]\tp
 \end{align*}
 As a normal subgroup of $G(\Phi,\Delta,R,\sigma)$, the last group is generated by the commutators $[x,y]$, where $x$,$y\in E(\Phi,\Delta,R,\sigma)\cup D^i(T(\Phi,R)G(\Delta,R))$. By item 2 of Lemma~\ref{ChevMats} all these commutators belong to
 $E(\Phi,\Delta,R,\sigma)D^{i+1}(T(\Phi,R)G(\Delta,R))$; hence it remains to prove that this subgroup is normal in $G(\Phi,\Delta,R,\sigma)$, but this follows immediately from the fact that the subgroup $E(\Phi,\Delta,R,\sigma)$ is normal in the group $G(\Phi,\Delta,R,\sigma)$, and $D^{i+1}(T(\Phi,R)G(\Delta,R))$  is normal in $T(\Phi,R)G(\Delta,R)$.
 
Further, note that
 $$
  E(\Phi,\Delta,R,\sigma)D^i(T(\Phi,R)G(\Delta,R))\le E(\Phi,\Delta,R,\sigma)D^{i-1}(G(\Delta,R))\tc
 $$
 and by Lemma~\ref{K1nilp} , if $i$ is sufficiently large, then the right-hand side is contained in
 $$
 E(\Phi,\Delta,R,\sigma)E(\Delta,R)=E(\Phi,\Delta,R,\sigma)\tp
 $$
 
 Now we consider the general case. From what has been proved previously, it follows
 that for $N_1$ sufficiently large, we have:
 $$
 \rho_J(D^{N_1}G(\Phi,\Delta,R,\sigma))\le E(\Phi,\Delta,\fact{R}{J},\fact{\sigma}{J})\tp
 $$  
 Since $E(\Phi,\Delta,R,\sigma)$ maps surjectively onto $E(\Phi,\Delta,\fact{R}{J},\fact{\sigma}{J})$, we obtain
$$
D^{N_1}G(\Phi,\Delta,R,\sigma)\le E(\Phi,\Delta,R,\sigma)G(\Phi,R,J)\tp
$$ 
Now, we prove that
 $$
 D^i( E(\Phi,\Delta,R,\sigma)G(\Phi,R,J))\le  E(\Phi,\Delta,R,\sigma)D^i(G(\Phi,R,J))\tp
 $$
 As before, it suffices to verify the fact that $E(\Phi,\Delta,R,\sigma)D^{i+1}(G(\Phi,R,J))$  is normal in $E(\Phi,\Delta,R,\sigma)G(\Phi,R,J)$ and contains the commutators of elements from $E(\Phi,\Delta,R,\sigma)\cup D^i(G(\Phi,R,J))$. Both claims follow from Lemma~\ref{standcomm} and the fact that by Lemma~\ref{ISquare} we have
 $$
 E(\Phi,R,J)\le E(\Phi,I^+\cap I^-)\le E(\Phi,\Delta,R,\sigma)\tp
 $$
 Using Lemmas~\ref{K1nilp} and \ref{ISquare} once again, we deduce that for a large $N_2$ 2 we have
\begin{align*}
 D^{N_2}( E(\Phi,\Delta,R,\sigma)G(\Phi,R,J))\le  E(\Phi,\Delta,R,\sigma)E(\Phi,R,J)\le E(\Phi,\Delta,R,\sigma)E(\Phi,I^+\cap I^-)=\\=E(\Phi,\Delta,R,\sigma)\tp
\end{align*}
 
 It remains to set $N=N_1+N_2$.
 \end{proof}
 
 {\lem\label{normality} The subgroup $E(\Phi,\Delta,R,\sigma)$ is normal in $G(\Phi,\Delta,R,\sigma)$.}
 \begin{proof}
 For a finitely generated ring, the statement follows from Lemma~\ref{DerivSer}, and every ring is an inductive limit of its finitely generated subrings. 
 \end{proof}
 
 {\prop\label{Transporter} We have the identity
 $$
 N_{G(\Phi,R)}(E(\Phi,\Delta,R,\sigma))=\Tran(E(\Phi,\Delta,R,\sigma),G(\Phi,\Delta,R,\sigma))\tp
 $$}
 \begin{proof}
 Obviously, the right-hand side is contained in the left-hand side. We prove the
 inverse inclusion. As in the previous lemma, without loss of generality we may assume
 that the ring $R$ is finitely generated.
 
 Then by Lemma~\ref{DerivSer}, we have:
 \begin{align*}
 \Tran(E(\Phi,\Delta,R,\sigma),G(\Phi,\Delta,R,\sigma))\le \Tran(D^N E(\Phi,\Delta,R,\sigma),D^N G(\Phi,\Delta,R,\sigma))=\\=\Tran(E(\Phi,\Delta,R,\sigma),E(\Phi,\Delta,R,\sigma))\tp
 \end{align*}
It remains to prove that if $g\in \Tran(E(\Phi,\Delta,R,\sigma),E(\Phi,\Delta,R,\sigma))$, then the same is true for $g^{-1}$.

First, if $g\in \Tran(E(\Phi,\Delta,R,\sigma),E(\Phi,\Delta,R,\sigma))$, then for every $k\in\N$ we have $g^k\in \Tran(E(\Phi,\Delta,R,\sigma),E(\Phi,\Delta,R,\sigma))$.

 Let $\tilde{g}$ be the operator on the space of matrices $M(\Lambda)$ given by:
\begin{align*}
\map{\tilde{g}}{&M(\Lambda)}{M(\Lambda)}\\
&m\mapsto gmg^{-1}\tp
\end{align*}
From the Cayley--Hamilton theorem, it follows that the operator $\tilde{g}^{-1}$ is a polynomial
of $\tilde{g}$. Set
$$
L=\{m\in M(\Lambda)\colon\forall\lambda\ne\lambda_0\gap m_{\lambda_0,\lambda}\in I^+ \,\&\, m_{\lambda,\lambda_0}\in I^-\}\tp
$$
 Since all operators $\tilde{g}^k$ preserve $E(\Phi,\Delta,R,\sigma)$, all polynomials of $\tilde{g}$ map $E(\Phi,\Delta,R,\sigma)$ to $L$. Since
$$
L\cap G(\Phi,R)=G(\Phi,\Delta,R,\sigma)\tc
$$
we obtain
$$
g^{-1}\in \Tran(E(\Phi,\Delta,R,\sigma),G(\Phi,\Delta,R,\sigma))\le \Tran(E(\Phi,\Delta,R,\sigma),E(\Phi,\Delta,R,\sigma))\tp
$$
 \end{proof}
 
The idea to use the nilpotent structure $K_1$ 1 for the computation of normalisers was
suggested by A. Stepanov in~\cite{StepStandard}. 
 
\section{The proof of Theorem 2 and its corollary}

\begin{proof}
 From Proposition~\ref{Transporter} it follows that $N_{G(\Phi,R)}(E(\Phi,\Delta,R,\sigma))$  consist exactly of elements $g\in G(\Phi,R)$ such that the group $E(\Phi,\Delta,R,\sigma)$  preserves the line spanned by the covector $g_{\lambda_0,*}$ modulo $I^+$, and preserves the line spanned by the vector $(g^{-1})_{*,\lambda_0}$ modulo $I^-$. Since the group $E(\Phi,\Delta,R,\sigma)$ is perfect, the preservation of lines implies
 the preservation of the vector and covector themselves. By writing out the meaning of
 this preservation for the generators of the group $E(\Phi,\Delta,R,\sigma)$ in terms of matrices, we obtain Theorem~\ref{Normalizer}. 
\end{proof}
 
To obtain the corollary we are interested in, we need the following two lemmas.
 
 {\lem\label{bilform} In the cases of the second type, for any $\lambda\in\Lambda$ the weight $-\lambda$ also belongs to $\Lambda$, and there exists $G(\Phi,R)$-invariant bilinear form on the module $V$ given by 
 $$
 h\left(\sum_{\lambda\in\Lambda} x_\lambda v^{\lambda},\sum_{\lambda\in\Lambda} y_\lambda v^{\lambda}\right)=\sum_{\lambda\in\Lambda}\pm x_\lambda y_{-\lambda}\tp
 $$
 }
 \begin{proof}
  The fact that $-\lambda_0\in \Lambda$ has already been mentioned. Thus, the first statement
  follows from the fact that $W(\Phi)$ acts transitively on weights. Next, consider the dual
  representation $V^*$ and its basis $\{(v^{\lambda})^*\}$ that is dual to the basis $\{v^\lambda\}$. Then the vector $(v^{\lambda})^*$ is a weight vector with the weight $-\lambda$. Therefore, the dual representation has the same set of weights; hence there is an isomorphism between $V$ and $V^*$ given by
 $$
 \map{\tilde{h}}{V}{V^*}\tc\qquad v^\lambda\mapsto \pm (v^{-\lambda})^*\tc
 $$
which is the same as an invariant bilinear form.
 \end{proof}
 
Note that for the case of $E_7$ such a form was written out explicitly in~\cite{LuzE7Inv}.  
 
 {\lem\label{qform} In the cases of the second type, there exists a quadratic form on the module~$V$:
 $$
 q\left(\sum_{\lambda\in\Lambda} x_\lambda v^{\lambda}\right)=\sum q_{\mu,\nu}x_{\mu}x_{\nu}\tc
 $$
 (where $(\mu,\nu)$ runs over the set of nonordered pairs of weights) that satisfies the following conditions:
 \begin{enumerate}
 \item $q(v)=0$ for any vector $v$ from the orbit of the highest weight vector {\rm (}i.e., for $v=g_{*,\lambda_0}$ for any $g\in G(\Phi,R)$; since the Weyl group acts transitively, the vectors $v=g_{*,\lambda}$ for the other weights $\lambda$ also belong to that orbit{\rm )};
 
 \item $q_{\lambda_0,-\lambda_0}=\pm 1$. 
 \end{enumerate}}
 
 In \cite{LuzEquations}, such equations on the orbit of the highest weight vector were called $\pi$-equations. 
 
 \begin{proof}
Let $\mu_1$,$\mu_2$,$\ldots$,$\mu_k$ be the shortest path in the weight graph such that $\mu_1\in \Lambda_2$ and $\mu_k=-\lambda_{0}$. Then $\gamma_i=\mu_i-\mu_{i+1}\in \Phi$ and $d(\mu_i,-\lambda_0)=k-i$.  
 
 We will prove by induction on $i$ that there exists a quadratic form
 $$
 q^{(i)}\left(\sum_{\lambda\in\Lambda} x_\lambda v^{\lambda}\right)=\sum q^{(i)}_{\mu,\nu}x_{\mu}x_{\nu}\tc
 $$
 where $(\mu,\nu)$ runs over the set of nonordered pairs of weights such that $d(\mu,-\lambda_0)\ge k-i$ and $d(\nu,-\lambda_0)\ge k-i$, that satisfies the following conditions:

 \begin{enumerate}
 \item $q^{(i)}(v)=0$for any vector $v$ from the orbit of the highest weight vector;
 
 \item $q^{(i)}_{\lambda_0,\mu_i}=\pm 1$. 
 \end{enumerate}  
 
The base of induction, $i=1$ follows from \cite{VavNumerology}. Theorem 2 of that paper implies that, since $d(\lambda_0,\mu_1)=2$, the set 
 $$
 \Omega(\lambda_0,\mu_1)=\{\nu\in\Lambda\colon d(\nu,\lambda_0)=d(\nu,\mu_1)=1\}\cup\{\lambda_0,\mu_1\}\sub \Lambda_1\cup\{\lambda_0,\mu_1\}
 $$
 is a {\it square}\footnote{A subset $\Omega\sub\Lambda$ is called a {\it square}, if $|\Omega|\ge 4$, and for any $\lambda\in\Omega$ exactly one of the differences $\{\lambda-\mu\colon \mu\in\Omega\sm\{\lambda\}\}$ is NOT a root.}. Hence the vectors from the orbit of the highest weight vector satisfy the corresponding square equation, which has the required shape. It only remains to observe that the weights from $\Lambda_1$ and the weight $\lambda_0$ are at a distance of at least $k$ from the weight $-\lambda_0$ in the weight graph. Indeed, otherwise the corresponding path goes through $\Lambda_2$ (by
 Lemma~\ref{alphaonecoefficient}, every root has coefficient of $\alpha^{(1)}$ not greater than one in absolute value; hence every vertex in the path is either in the same component as the previous one, or in the adjacent component; in other words, the path cannot jump through $\Lambda_2$), and we get a shorter path to $\Lambda_2$.
 
Now we do the induction step from $i$ to $i+1$. Set
 $$
 q^{(i+1)}(v)=q^{(i)}(x_{\gamma_i}(1)v)\tp
 $$
 After this transformation, the minimal distance from a weight participating in the formula
 to the weight $-\lambda_0$ decreases at most by 1. Obviously, the new form still annihilates the vectors from the orbit of the highest weight vector. Since $d(\mu_{i+1},-\lambda_0)=k-i-1$, the weight $\mu_{i+1}$ does not participate in the formula for $q^{(i)}$; hence the monomial $x_{\lambda_0}x_{\mu_i}$ of the form $q^{(i)}$s the only one that contributes to the coefficient of $x_{\lambda_0}x_{\mu_{i+1}}$ in the form $q^{(i+1)}$. Therefore, this coefficient is equal to $\pm 1$.
 \end{proof}

{\cor\label{DerivSer2} (to Theorem $\ref{Normalizer}$) Let $R$ be a finitely generated ring, and let $\sigma=(I^+,I^-)$ be a pair of ideals of $R$. Then there exists a natural number $N$ such that
  $$
 D^N N_{G(\Phi,R)}\left(E(\Phi,\Delta,R,\sigma)\right)=E(\Phi,\Delta,R,\sigma)\tp
 $$}
 \begin{proof}
 In the cases of the first type, the statement follows from Theorem~\ref{Normalizer} and Lemma~\ref{DerivSer}. We consider the case of the second type.
 
 By Lemma~\ref{DerivSer}, it suffices to prove the inclusion
 $$
 \left[N_{G(\Phi,R)}\left(E(\Phi,\Delta,R,\sigma)\right),N_{G(\Phi,R)}\left(E(\Phi,\Delta,R,\sigma)\right)\right]\le G(\Phi,\Delta,R,\sigma)\tp
 $$ 
We prove the inclusion
 $$
 \left[N_{G(\Phi,R)}\left(E(\Phi,\Delta,R,\sigma)\right),N_{G(\Phi,R)}\left(E(\Phi,\Delta,R,\sigma)\right)\right]\le G(\Phi,\Delta,R,(R,I^-))\tc
 $$ 
 The proof of the inclusion into $G(\Phi,\Delta,R,(I^+,R))$ is similar (the form $q$ from Lemma~\ref{qform} should be transferred to the dual representation via the isomorphism described in Lemma~\ref{bilform}).   
 
 Replacing the ring by its quotient $\fact{R}{I^-}$, we may assume that $I^-=(0)$. Therefore,
 we need to prove the inclusion
 $$
 \left[N_{G(\Phi,R)}\left(E(\Phi,\Delta,R,\sigma)\right),N_{G(\Phi,R)}\left(E(\Phi,\Delta,R,\sigma)\right)\right]\le P\tp
 $$ 
 
 Since the ring $R$ is finitely generated, it is Noetherian; hence it is a direct product of rings with connected spectra. Thus, without loss of generality, we may assume that the spectrum of $R$ is connected.
 
 {\bf Case 1}. $I^+\ne (0)$.
 
 We claim that in this case, even a stronger inclusion holds:
 $$
 N_{G(\Phi,R)}\left(E(\Phi,\Delta,R,\sigma)\right)\le P\tp
 $$ 
 Indeed, let
 $g\in N_{G(\Phi,R)}\left(E(\Phi,\Delta,R,\sigma)\right)$. Then $g^{-1}\in N_{G(\Phi,R)}\left(E(\Phi,\Delta,R,\sigma)\right)$, and by Theorem~\ref{Normalizer} we have
 $$
 g_{\lambda,\lambda_0}\in I^-=(0)\quad\forall\lambda\in \Lambda\sm\{\lambda_0,-\lambda_{0}\}\tp
 $$
 Next, consider the form $q$ from Lemma~\ref{qform}. Since it annihilates the vectors from the
 orbit of the highest weight vector, the coefficients $q_{\lambda,\lambda}$($=q(v^\lambda)$) vanish; hence, using that $g_{*,\lambda_0}=g_{\lambda,\lambda_0}v^{\lambda_0}+g_{-\lambda_0,\lambda_0}v^{-\lambda_0}$,  we obtain
 $$
 g_{\lambda,\lambda_0}g_{-\lambda_0,\lambda_0}=\pm q(g_{*,\lambda_0})=0\tp
 $$ 
 The pair $g_{\lambda_0,\lambda_0}$,$g_{-\lambda_0,\lambda_0}$ is unimodular. Hence every prime ideal of the ring $R$ contains exactly one of the elements $g_{\lambda_0,\lambda_0}$ and $g_{-\lambda_0,\lambda_0}$. Since the spectrum of the ring $R$ is
 connected, one of those elements is invertible, and the other (since $g_{\lambda_0,\lambda_0}g_{-\lambda_0,\lambda_0}=0$) is
 equal to zero. By Theorem~\ref{Normalizer} we have the inclusion
 $$
 g_{-\lambda_0,\lambda_0}I^+\sub I^-=0\tp
 $$
 Hence the element $g_{-\lambda_0,\lambda_0}$ is not invertible, i.e., it is equal to zero. Therefore, $g\in P$.
 
 {\bf Case 2}. $I^+= (0)$.
 
 First we assume that $g\in N_{G(\Phi,R)}\left(E(\Phi,\Delta,R,\sigma)\right)$. Then $g^{-1}\in N_{G(\Phi,R)}\left(E(\Phi,\Delta,R,\sigma)\right)$, and by Theorem~\ref{Normalizer} we have:
\begin{align*}
 g_{\lambda,\lambda_0}\in I^-=(0)\quad\forall\lambda\in \Lambda\sm\{\lambda_0,-\lambda_{0}\}\tc\\
  g_{\lambda_0,\lambda}\in I^+=(0)\quad\forall\lambda\in \Lambda\sm\{\lambda_0,-\lambda_{0}\}\tp
\end{align*}

The same is true for $g^{-1}$; hence for $\lambda\in \Lambda\sm\{\lambda_0,-\lambda_{0}\}$, we have
$$
g_{-\lambda_0,\lambda}=h(g_{*,\lambda},v^{\lambda_0})=h(gv^{\lambda},v^{\lambda_0})=h( v^{\lambda},g^{-1}v^{\lambda_0})=(g^{-1})_{-\lambda,\lambda_0}=0\tc
$$
where $h$ is a form occurring in Lemma~\ref{bilform}. Similarly, we deduce that
$
g_{\lambda,-\lambda_0}=0\tp
$
Therefore, the submodule generated by the vectors $v^{\lambda_0}$ and $v^{-\lambda_0}$, is invariant with respect to the action of the group $N_{G(\Phi,R)}\left(E(\Phi,\Delta,R,\sigma)\right)$, and the element $g$ acts on it by the
following matrix:
$$
\begin{pmatrix}
g_{\lambda_0,\lambda_0} & g_{\lambda_0,-\lambda_0}\\
g_{-\lambda_0,\lambda_0} & g_{-\lambda_0,-\lambda_0}
\end{pmatrix}\tp
$$
Arguing as in the first case, we see that each column and each row of this matrix contain
an invertible entry and a zero entry. Therefore, this matrix can be of one of the following
shapes:
$$
\begin{pmatrix}
g_{\lambda_0,\lambda_0} &           0                 \\
           0            & g_{-\lambda_0,-\lambda_0}
\end{pmatrix}\text{or}
\begin{pmatrix}
           0             & g_{\lambda_0,-\lambda_0}\\
g_{-\lambda_0,\lambda_0} &          0            
\end{pmatrix}\tp
$$ 
If
 $$
 g\in \left[N_{G(\Phi,R)}\left(E(\Phi,\Delta,R,\sigma)\right),N_{G(\Phi,R)}\left(E(\Phi,\Delta,R,\sigma)\right)\right]\tc
 $$
  then the matrix has the first shape, i.e. $g_{-\lambda_0,\lambda_0}=0$ and $g\in P$.
\end{proof}  
 
 \section{Root type elements}
 
Let $\Root(-)$ be the smallest closed subscheme (over $\Z$) in $G(\Phi,-)$, such that for any $R$ and any $h\in G(\Phi,R)$ we have
$$
x_{\alpha_1}(1)^h\in \Root(R)\tp
$$ 

Obviously, one can replace $\alpha_1$ with any other root because all elements $x_{\alpha}(1)$ are conjugate by elements of the Weyl group.

The elements of $\Root(R)$ are what we call the {\it root type elements}.

Let $\Z[G]$ be the ring of regular functions on the scheme $G(\Phi,-)$, and let $\Z[\Root]$ be the ring of regular functions on the scheme $\Root(-)$. By definition, the second ring is a quotient of the first one. However, using the following lemma, we will also view the ring $\Z[\Root]$ as a subring of $\Z[G]$. 

We denote by $g_{\gen}\in G(\Phi,\Z[G])$ and $r_{\gen}\in \Root(\Z[\Root])$ the generic elements of the corresponding schemes.

{\lem\label{ZRootToZG} There exists an injective map
$
\map{i}{Z[\Root]}{\Z[G]}
$
such that 
$
i(r_{\gen})=x_{\alpha_1}(1)^{g_{\gen}}\tp
$}

\begin{proof}
The collection of maps
\begin{align*}
G(\Phi,R)\to \Root(R)\\
g\mapsto x_{\alpha_1}(1)^g
\end{align*}
gives a morphism of schemes $G(\Phi,-)\to \Root(-)$. Take $i$ to be the corresponding homomorphism of rings of regular functions. The relation $i(r_{\gen})=x_{\alpha_1}(1)^{g_{\gen}}$ is true in this
case by definition. It remains to prove that the homomorphism $i$ is injective.

Let $I=\Ker i$, and let $\tilde{\Root}(-)$ be the subscheme of $\Root(-)$ determined by the ideal $I$. Then by the definition of $I$  we have
$
x_{\alpha_1}(1)^{g_{\gen}}\in \tilde{\Root}(\Z[G])\tp
$
Hence by the universal property of $g_{gen}$, we have 
$
x_{\alpha_1}(1)^{h}\in \tilde{\Root}(R)\tp
$
for any ring $R$ and any $h\in G(\Phi,R)$.  However, $\Root(-)$ ) is the
smallest subscheme with this property. Hence  $\tilde{\Root}(-)$ coincides with $\Root(-)$, i.e. $I=0$.
\end{proof}

We provide several important examples of root type elements.
{\lem\label{examples} Let $R$ be a commutative ring; then the following holds
\begin{enumerate}
\item $x_\alpha(t)\in \Root(R)$ for any $\alpha\in\Phi$, $t\in R$. In particular, the identity element is a root type element.
\item An element conjugate to a root type element is a root type element.
\item Let $\alpha$,$\beta\in \Phi$ be such that $\angle(\alpha,\beta)=\tfrac{\pi}{3}$. Then 
$
x_\alpha(\xi)x_{\beta}(\zeta)\in \Root(R)\tp
$
\end{enumerate}
}

\begin{proof}
\begin{enumerate}
\item  It suffices to consider the case where $t$ is a free variable ($R=\Z[t]$). We
need to check that the element $x_\alpha(t)$ satisfies certain equations. Since the ring $\Z[t]$ is embeded into the ring $\Z[t,t^{-1}]$, we can replace one with another. In the
group $G(\Phi,\Z[t,t^{-1}])$ the element $x_\alpha(t)$ is conjugate to the element $x_{\alpha_1}(1)$; hence it satisfies those equations.

\item  It suffices to consider the element
$
(r_{\gen})^{g_{\gen}}\in G(\Phi,\Z[\Root]\otimes_{\Z}\Z[G])\tp
$
As before, it suffices to map our ring to another ring injectively with the image of the element $r_{\gen}$ (and hence of the element $(r_{\gen})^{g_{\gen}}$) being conjugate to $x_{\alpha_1}(1)$. Those requirements are satisfied if we take the map
$$
\map{i\otimes\Z[G]}{\Z[\Root]\otimes_{\Z}\Z[G]}{\Z[G]\otimes_{\Z}\Z[G]}\tc
$$ 
where $i$ is from Lemma~\ref{ZRootToZG}. This map is injective because $i$ is injective, and the
scheme $G(\Phi,-)$ is flat, i.e $\Z[G]$ is a flat $\Z$-module.

\item As before, it suffices to consider the ring $R=\Z[\xi,\zeta,\zeta^{-1}]$.  Note that in the group $\SL(3,R)$) we have
$$
\begin{pmatrix}
1 & 0 & 0\\
0 & 1 & 0\\
0 & \xi\zeta^{-1} & 1
\end{pmatrix}
\begin{pmatrix}
1 & \xi & \zeta\\
0 & 1 & 0\\
0 & 0 & 1
\end{pmatrix}
\begin{pmatrix}
1 & 0 & 0\\
0 & 1 & 0\\
0 & -\xi\zeta^{-1} & 1
\end{pmatrix}=\begin{pmatrix}
1 & 0 & \zeta\\
0 & 1 & 0\\
0 & 0 & 1
\end{pmatrix}\tp
$$
Applying the corresponding map from $\SL(3,R)$ to $G(\Phi,R)$, we see that the
element $x_\alpha(\xi)x_\beta(\zeta)$ is conjugate to the element $x_\beta(\zeta)$; hence it is conjugate to $x_{\alpha_1}(1)$.
\end{enumerate}
\end{proof}

{\lem\label{RootTypeAsExponent}
Every root type element is an exponential of a certain element of the Lie algebra of the group $G(\Phi,R)$, Moreover, in terms of matrices in the representation $V$ taking the exponential is addition of the identity matrix.}
\begin{proof}
It suffices to prove the statement for the element $r_{\gen}\in \Root(\Z[\Root])$ .We view $\Z[\Root]$ as a subring in $\Z[G]$ (Lemma~\ref{ZRootToZG}). Over $\Z[G]$, the element $r_{\gen}$ becomes conjugate to the element $x_{\alpha_1}(1)$; hence the statement is true for it (the second part holds true because the representation is minuscule). It remains to note that after subtraction of the identity matrix, all the entries remain in the ring $\Z[\Root]$; hence the corresponding element of the Lie algebra is defined over $\Z[\Root]$.
\end{proof}

{\lem\label{rootequation} Let $g\in\Root(R)$. Then the following holds true:
\begin{enumerate}
\item $g_{\lambda,\mu}=0$ for any pair of weights $\lambda,\mu$ such that $d(\lambda,\mu)\ge 2$.

\item $g_{\lambda,\mu}=\pm g_{\rho,\sigma}$ for any $\lambda$, $\mu$, $\rho$, $\sigma\in\Lambda$ such that $d(\lambda,\mu)=d(\rho,\sigma)=1$ and $\lambda-\mu=\rho-\sigma$.
\end{enumerate}
}
\begin{proof}
The two statements follow from Lemma~\ref{RootTypeAsExponent} and the Cartan decomposition.
\end{proof}

 \section{Extraction of elementary root elements}
 
We need the following obvious lemma.
 
 {\lem\label{expand} Let $\lambda\in\Lambda$ and $\alpha$,$\beta_1$,$\ldots$,$\beta_k\in\Phi$ be such that $\lambda+\alpha\notin\Lambda$, and $\lambda-\beta_i\in\Lambda$ for all~$i$. Then for any $\zeta$,$\xi_1$,$\xi_k\in R$ we have
 $$
 \left[\prod_{i=1}^k x_{\beta_i}(\xi_i),x_\alpha(\zeta)\right]=\prod_{i=1}^k[x_{\beta_i}(\xi_i),x_\alpha(\zeta)]\tp
 $$}
 
 \subsection{Extraction from the subgroups $P$ and $P^-$}
 
The proposition below follows from the paper \cite{StavNormStr}. However, in our case the proof is quite simple, so we give it here.
 
 {\prop\label{PExtraction}
Let $H$ H be an overgroup for $E(\Delta,R)$.
\begin{enumerate}
\item
 Let $I^+\unlhd R$ be such that
 $
 (H\cap P)\sm G(\Phi,\Delta,R,(I^+,0))\ne \emp\tp
 $
 Then $H$ contains an element $x_\alpha(\xi)$, where $\alpha\in \Omega^+$ and $\xi\in R\sm I^+$.
 \item The same is true for $P^-$ and $\Omega^-$. 
 \end{enumerate} }
 \begin{proof}
 
We prove Item 1, the proof of the second item is similar. 
 
Suppose the contrary; then the first component of $\lev(H)$ is contained in $I^+$. Without
loss of generality, we may assume that it is equal to $I^+$ (if we replace $I^+$ with this
component, then the assumption remains true).
 
Let
$
g\in (H\cap P)\sm G(\Phi,\Delta,R,(I^+,0))\tc
$
then $g=g_1g_2$, where $g_1\in U$ and $g_2\in L$. Let
$$
g_1=\prod_{i=1}^l x_{\beta_i'}(\xi_i)\tc\qquad \text{where } \beta_i'\in \Omega^+\tp
$$

Since $g\notin G(\Phi,\Delta,R,(I^+,0))$, without loss of generality we may assume that $\xi_1\notin I^+$. Take $\alpha\in\Delta$ such that $\beta_1'+\alpha\in\Phi$ (Lemma~\ref{RootInDeltaComb}). Then we have:
$$
[g_1,x_\alpha(1)]=[gg_2^{-1},x_\alpha(1)]=\phan^g[g_2^{-1},x_\alpha(1)]\cdot [g,x_\alpha(1)]\in H\tc
$$
because $g$, $x_\alpha(1)\in H$, and
$$
[g_2^{-1},x_\alpha(1)]\in [L,E(\Delta,R)]=[G(\Phi,\Delta,R,(0,0)),E(\Phi,\Delta,R,(0,0))]=E(\Phi,\Delta,R,(0,0))\le H\tp
$$
On the other hand, by Lemma~\ref{expand} we get:
$$
[g_1,x_\alpha(1)]=\prod_{i=1}^l[x_{\beta_i'}(\xi_i),x_\alpha(1)]\tp
$$
Without loss of generality, we may assume that the first $k$ commutators and only they are nontrivial, i.e.,
$$
[g_1,x_\alpha(1)]=\prod_{i=1}^kx_{\beta_i}(\pm\xi_i)\tc\qquad \text{where } \beta_i=\beta'_i+\alpha\tp
$$

Therefore, $H$ contains an element of the form
$$
\prod_{i=1}^kx_{\beta_i}(\xi_i)\tc
$$
where $\beta_i\in\Omega^+$ and $\xi_1\notin I^+$. Among all such elements of $H$ consider those that have the smallest $k$, and among them, take the one with the greatest sum of the $\beta_i$. Then none of the $\xi_i$ belongs to $I^+$ (otherwise, the corresponding root elements belong to $H$, and we can remove them to reduce $k$), and all the $\beta_i$ are distinct (otherwise we reduce $k$, by using additivity). If at least one of the $\beta_i$ is not the maximal root, then we can take a simple root $\alpha_m$ such that $\beta_i+\alpha_m\in\Phi$ ($\alpha_m\ne \alpha^{(1)}$ because $\beta_i\in\Omega^+$; hence $\alpha_m\in\Delta$).  In this case, we see that the element
$$
\left[\prod_{i=1}^kx_{\beta_i}(\xi_1),x_{\alpha_m}(1)\right]\in H
$$
has the same form by Lemma~\ref{expand}. Furthermore, we either reduce $k$, or retain $k$ the same, but increase the sum of $\beta_i$, which contradicts our choice

Therefore, our element is $x_\delta(\xi_1)$, where $\delta$ is the maximal root, and $\xi_1\notin I^+$.
\end{proof}
 
  \subsection{Extraction from the subgroups $P_\lambda$ and $P^-_\lambda$}
  {\lem\label{UnipExtraction} Let $H$ be an overgroup of $E(\Delta,R)$, of level $\sigma=(I^+,I^-)$, and let $\lambda_1\in\Lambda_1$. Then the following inclusion holds:
  $$
  H\cap U_{\lambda_1}\sub G(\Phi,\Delta,R, \sigma)\tp
  $$
  } 
  \begin{proof}
  Let $g\in H\cap U_{\lambda_1}\sm G(\Phi,\Delta,R, \sigma)$. 
  
  {\bf Case 1}. $g\notin G(\Phi,\Delta,R, (I^+,R))$.
    
In this case, we have
  $$
  g=\prod_{i=0}^l x_{\beta_i'}(\xi_i)\tc\qquad \text{where } \beta_i'\in \Sigma_{\lambda_1}\tp
  $$
  We may assume that $\beta_0'=\lambda_1-\lambda_0$, and $x_{\beta_1'}(\xi_1)\notin  G(\Phi,\Delta,R, (I^+,R))$.  Then$\beta_1'\in \Sigma_{\lambda_1}^+$ and $\xi_1\notin I^+$.
  
  By Item 3 of Lemma~\ref{SigmaStructure} there exists $\gamma\in \Delta\cap\Delta_{\lambda_1}$ such that $\beta_1+\gamma\in\Phi$  (hence, in fact, $\beta_1+\gamma\in\Sigma_{\lambda_1}^+$), but $\gamma+\beta_0\notin \Phi$. Then
  $
  [g,x_{\gamma}(1)]\in H\tc
  $
  On the other hand, by Lemma~\ref{expand} we have:
  $$
  [g,x_{\gamma}(1)]=\prod_{i=0}^l[x_{\beta_i'}(\xi_i),x_{\gamma}(1)]=\prod_{i=1}^l[x_{\beta_i'}(\xi_i),x_{\gamma}(1)]\tp
  $$
  Without loss of generality, we may assume that the first $k$ commutators and only they
  are nontrivial, i.e.,
  $$
 [g,x_{\gamma}(1)]=\prod_{i=1}^k x_{\beta_i}(\pm\xi_i)\tc
  $$
  where all the $\beta_i$ belong to $\Sigma_{\lambda_1}^0\sqcup \Sigma_{\lambda_1}^+$ (because by Item 1 of Lemma~\ref{SigmaStructure} $\beta'_0$ is the only element of $\Sigma_{\lambda_1}^-$); next, $\beta_1\in \Sigma_{\lambda_1}^+$ and $\xi_1\notin I^+$. Arguing as in the end of the proof of Proposition~\ref{PExtraction}, we get a contradiction.
  
  {\bf Case 2}. $g\in G(\Phi,\Delta,R, (I^+,R))$.
  
  Similarly,
  $$
  g=\prod_{i=0}^l x_{\beta_i'}(\xi_i)\tp
  $$
  However, now we have $\xi_i\in I^+$ ($1\le i\le l)$, but $\xi_0\notin I^-$. Thus, all the factors except the first one belong to $H$; hence we deduce that 
  $
  x_{\beta_0'}(\xi_0)\in H\tc
  $
 which contradicts the fact that $\xi\notin I^-$.
  \end{proof}
 
{\lem \label{MagicManipulation}  Let $I^+\unlhd R$, $\lambda_1\in\Lambda_1$, and let
 $$
 g\in(P_{\lambda_1}\cap\Root(R))\sm G(\Phi,\Delta,R,(I^+,R))\tc
 $$
 
 Then there exists $\gamma\in(\Delta\cap \Delta_{\lambda_1})'$ {\rm (}see Item (4) of Lemma~$\ref{SigmaStructure}${\rm )} such that
 $$
 gx_{\gamma}(1)g^{-1}\notin G(\Phi,\Delta,R,(I^+,R))\tp
 $$}
 \begin{proof}
 
 First, we prove the following statement:
 $$
 gE((\Delta\cap \Delta_{\lambda_1})',R)g^{-1}\not\sub G(\Phi,\Delta,R,(I^+,R))\tp
 $$
 
 Suppose the contrary; then our assumption means exactly that the group \linebreak$E((\Delta\cap \Delta_{\lambda_1})',R)$ stabilizes the line spanned by the covector $g_{\lambda_0*}$, modulo $I^+$. Since this group is perfect (because $(\Delta\cap \Delta_{\lambda_1})'$ is an irreducible system distinct from $A_1$), it must stabilize
 the covector itself.
  
Further, since $g\notin G(\Phi,\Delta,R,(I^+,R))$, there exists $\mu\in\Lambda\sm\{\lambda_0\}$ such that $g_{\lambda_0,\mu}\notin I^+$. By Lemma~\ref{rootequation} $\mu\in\Lambda_1$, and since $g\in P_{\lambda_1}$, we see that $\mu\ne \lambda_1$. Take $\nu$ from Item 4 of Lemma~\ref{SigmaStructure}, and set
  $
  \tilde{\gamma}=\mu-\nu\in (\Delta\cap\Delta_{\lambda_1})'\tp
  $

Since the element $g_{\lambda_0*}$ stabilizes the covector $x_{\tilde{\gamma}}(1)$ we obtain the following congruence
  $$
  (gx_{\tilde{\gamma}}(1))_{\lambda_0,\nu}\equiv g_{\lambda_0,\nu}\mod I^+\tp
  $$
 However, we also have the identity
  $$
  (gx_{\tilde{\gamma}}(1))_{\lambda_0,\nu}=g_{\lambda_0,\nu}\pm g_{\lambda_0,\mu}\tp
  $$
  Hence $g_{\lambda_0,\mu}\in I^+$, which contradict the choice of $\mu$.
  
  Therefore, for some $\gamma\in(\Delta\cap \Delta_{\lambda_1})'$ and some $t\in R$, we have
  $$
 gx_{\gamma}(t)g^{-1}\notin G(\Phi,\Delta,R,(I^+,R))\tp
 $$
 We also know that
 $$
 gx_{\gamma}(t)g^{-1}=e+t(gx_{\gamma}(1)g^{-1}-e)\tc
 $$
Hence we obtain
 $$
 gx_{\gamma}(1)g^{-1}\notin G(\Phi,\Delta,R,(I^+,R))\tp
 $$
\end{proof}  
 
{\lem\label{GotoU} Let $\sigma=(I^+,I^-)$ be the pair of ideals of the ring $R$, let $\lambda_1\in\Lambda_1$, and let
$$
g\in (L_{\lambda_1}\cap \Root(R))\sm G(\Phi,\Delta,R,\sigma)\tp
$$
Then there exists $\gamma\in \Sigma_{\lambda_1}^0$ such that
$
 g^{-1}x_{\gamma}(1)g\notin G(\Phi,\Delta,R,\sigma)\tp
$
}
\begin{proof}
{\bf Case 1}. $g\notin G(\Phi,\Delta,R, (R,I^-))$.

In this case, there exists a weight $\mu$ such that
$
g_{\mu,\lambda_0}\notin I^-\tp
$
By Item 1 of Lemma~\ref{rootequation}, $\mu\in\Lambda_1$, and, since $g\in L_{\lambda_1}$, we see that $\mu\ne\lambda_1$. Next, let $w\in W(\Phi)$ be the reflection with respect to $\lambda_0-\lambda_1$. Note that the roots $\lambda_0-\mu$ and $\lambda_0-\lambda_1$ are not orthogonal. Indeed, otherwise we obtain:
$$
\lambda_0-\mu=w(\lambda_0-\mu)=\lambda_1-w(\mu)\tp
$$ 
Then, applying Item 2 of Lemma~\ref{rootequation} and the fact that $g\in L_{\lambda_1}$, we deduce that
$
g_{\mu,\lambda_0}=g_{w(\mu),\lambda_1}=0\tp
$
 
Therefore, their inner product should be equal to one (it cannot be equal to minus
one because they both belong to $\Omega^+$), and we can set
$$
\gamma=\lambda_1-\mu=(\lambda_0-\mu)-(\lambda_0-\lambda_1)\in\Phi\tp
$$
By definition, $\gamma\in\Sigma_{\lambda_1}$, and since $\mu\in\Lambda_1$, we have $\gamma\in\Sigma_{\lambda_1}^0$. We show that it satisfies the
requirement.

Indeed, the column $(x_\gamma(1)g)_{*,\lambda_0}$ is not a multiple of the column $g_{*,\lambda_0}$ modulo $I^-$, because $g_{\lambda_1,\lambda_0}=0$, but $(x_\gamma(1)g)_{\lambda_1,\lambda_0}\notin I^-$. This means exactly that
 $$
 g^{-1}x_\gamma(1)g\notin G(\Phi,\Delta,R,(R,I^-))\tp
 $$
 
 {\bf Case 2}. $g\notin G(\Phi,\Delta,R, (I^+,R))$.

In this case $g^{-1}\notin G(\Phi,\Delta,R, (I^+,R))$, i.e., the set
$$
M=\{\mu\in \Lambda_1 \colon g^{-1}_{\lambda_0,\mu}\notin I^+\}
$$
is nonempty.

Consider the strict partial order $\rightarrow$ on the set $M$ given by
$$
\mu_1\rightarrow \mu_2\Equal \exists\ovl{\eps}\in\fact{R}{I^+}\colon \left(\ovl{g^{-1}_{\lambda_0,\mu_1}}=\ovl{\eps}\ovl{g^{-1}_{\lambda_0,\mu_2}}\gap \text{and} \gap \ovl{\eps} \ovl{g^{-1}_{\lambda_0,\mu_1}}=0 \right)\tc
$$
where the bar stand for the image in $\fact{R}{I^+}$.

Multiplying the corresponding $\ovl{\eps}$, we get transitivity; antireflexivity follows by the
definition of $M$.

Let $\mu\in M$ be a maximal element (i.e., such that there are no arrows from it). As
in the first case, we prove that the root $\alpha=\lambda_0-\mu$has the inner product with $\lambda_0-\lambda_1$ equal to one. Take $\gamma$ from item (5) of Lemma~\ref{SigmaStructure}. Then $\alpha+\gamma\in\Omega^+$, and by Lemma~\ref{OmegaAndHighestWeight} we have:
$$
\nu=\mu-\gamma=\lambda_0-(\alpha+\gamma)\in\Lambda\tp
$$
Clearly, $\nu\in \Lambda_1$.

Assume that
$
 g^{-1}x_{\gamma}(1)g\in G(\Phi,\Delta,R,\sigma)\tc
$
i.e., the element $x_{\gamma}(1)$ stabilises the line
spanned by the covector $g^{-1}_{\lambda_0,*}$ modulo $I^+$. Thus, it multiplies this covector by a scalar, which we denote by$(1+\ovl{\eps})\in \fact{R}{I^+}$. 

We show that $\nu\in M$ and $\mu\rightarrow \nu$. If we do this, then we get a contradiction.

First, we have:
$$
g^{-1}_{\lambda_0,\nu}\pm g^{-1}_{\lambda_0,\mu}=(g^{-1}_{\lambda_0,*}x_{\gamma}(1))_{\nu}\equiv(1+\eps)g^{-1}_{\lambda_0,\nu}\mod I^+\tc
$$ 
i.e.,
$
\ovl{g^{-1}_{\lambda_0,\mu}}=(\pm\ovl{\eps})\ovl{g^{-1}_{\lambda_0,\nu}}\tp
$
In particular, this implies that $\nu\in M$. 

Second, we have:
$$
g^{-1}_{\lambda_0,\mu}=(g^{-1}_{\lambda_0,*}x_{\gamma}(1))_{\mu}\equiv (1+\eps)g^{-1}_{\lambda_0,\mu}\tc
$$
i.e.,
$
\ovl{\eps} \ovl{g^{-1}_{\lambda_0,\mu}}=0\tp
$
The lemma is proved.
\end{proof} 
 
 {\prop\label{PlambdaExtraction}
 Let $H$ be an overgroup of $E(\Delta,R)$.
\begin{enumerate}
\item
  Let $I^+\unlhd R$, and $\lambda_1\in\Lambda_1$ be such that there exists an element
 $$
 g\in(H\cap P_{\lambda_1}\cap\Root(R))\sm G(\Phi,\Delta,R,(I^+,R))\tp
 $$ 
 Then $H$ contains an element of the form $x_\alpha(\xi)$, where $\alpha\in \Omega^+$ and $\xi\in R\sm I^+$.
 \item The same is true for $P^-_{\lambda_1}$ and $\Omega^-$. 
 \end{enumerate} }
 \begin{proof}
We prove Item (1); the proof of Item (2) is similar.
 
Suppose the contrary; then the first component of $\lev(H)$) is contained in $I^+$. Without
loss of generality, we may assume that it is equal to $I^+$ (if we replace $I^+$ with this
component, then the assumption remains true). Let $I^-$ be the second component of $\lev(H)$ and let $\sigma=(I^+,I^-)$.

By Lemma~\ref{MagicManipulation} the element
$
g_1=gx_{\gamma_1}(1)g^{-1}
$
satisfies the same assumptions as the element $g$ for a suitable $\gamma_1\in(\Delta\cap\Delta_{\lambda_1})'$.

Let $g=ul$ and $g_1=u_1l_1$,
where $l$, $l_1\in L_{\lambda_1}$ and $u$, $u_1\in U_{\lambda_1}$. Then, taking the
projection to $L_{\lambda_1}$ in the definition of $g_1$, we obtain
$$
l_1=l^{-1}x_{\gamma_1}(1)l\in\Root(R)\tp
$$

{\bf Case 1}. $l_1\notin G(\Phi,\Delta,R,\sigma)$.

In this case, applying Lemma~\ref{GotoU}, we obtain the element
$$
h=l_1^{-1}x_{\gamma_2}(1)l_1\in U_{\lambda_1}\sm G(\Phi,\Delta,R,\sigma)\tp
$$
Next, since the group $U_{\lambda_1}$ is Abelian, we obtain
$
h=g_1^{-1}x_{\gamma_2}(1)g_1\in H\tp
$ 
Applying Lemma~\ref{UnipExtraction}, we get a contradiction.

{\bf Case 2}. $l_1\in G(\Phi,\Delta,R,\sigma)$.

Applying Lemma~\ref{MagicManipulation} once again, we obtain the element
$
g_2=g_1x_{\gamma_2}(1)g_1^{-1}
$
that satisfies the same assumptions as the initial element $g$. Let $g_2=u_2l_2$. As for the element $g_1$, we see that
$
l_2=l_1x_{\gamma_2}(1)l_1^{-1}\tp
$
By assumption, we have
$
l_1\in G(\Phi,\Delta,R,\sigma)
$;
hence by Lemma~\ref{normality}, we have
$$
l_2\in E(\Phi,\Delta,R,\sigma)\sub H\cap G(\Phi,\Delta,R,\sigma)\tp
$$
Hence $u_2\in H\sm G(\Phi,\Delta,R,\sigma)$, and we can apply Lemma~\ref{UnipExtraction} again.
 \end{proof}
 
\subsection{Extraction from a congruence subgroup of a nilpotent level}

{\lem Let $S$ be a commutative ring, and let $I\unlhd S$ be an ideal such that the following holds true.
\begin{enumerate}
\item $I^2=0$

\item The idea $I$ is finitely generated.

\item As an Abelian group, $S=\Z\+ I$, where $(1,0)$ is the identity element of $S$.

\item The ideal $I$ is torsion free as an Abelian group.
\end{enumerate}

Further, let $\xi\in I\sm\{0\}$. Then there exists a ring homomorphism
$$
\map{\ph}{S}{\fact{\C[\eps]}{(\eps^2)}}
$$
such that $\ph(\xi)\ne 0$.}

\begin{proof}
Since the ideal $I$ acts on itself by zero and is finitely generated as an ideal, it is
finitely generated as an Abelian group. Hence, as an Abelian group,
$$
I\simeq \bigoplus_{i=1}^N \Z\tp 
$$

Clearly, any subgroup of $I$ is an ideal. Choosing a direct summand with the projection
of $\xi$ to it being nonzero, and taking the quotient of $S$ by the sum of all other summands,
we may assume that $N=1$. In this case
$$
S\simeq \fact{\Z[\eps]}{(\eps^2)}\tc
$$
which can naturally be embedded into $\fact{\C[\eps]}{(\eps^2)}$.
\end{proof}

{\lem Let $R$ be an arbitrary commutative ring, and let $\B\unlhd R$ be an ideal such that $\B^2=0$. Next, let
$$
g\in G(\Phi,R,\B)\tp
$$
Then $g_{\lambda_0,\mu_0}=0$ for any pair of weights $\lambda_0,\mu_0$ such that $d(\lambda_0,\mu_0)\ge 2$.}

\begin{proof}
 Consider the ring
$$
\tilde{S}=\fact{\Z[\{a_{\lambda,\mu}\}_{\lambda,\mu\in\Lambda}]}{\<\{a_{\lambda,\mu}a_{\lambda',\mu'}\colon\lambda,\mu,\lambda',\mu'\in \Lambda\}\>}\tc
$$
and the ideal $\tilde{I}\unlhd \tilde{S}$ generated by all the $a_{\lambda\mu}$. Clearly, they satisfy the assumptions in the previous lemma.

Set
$
S=\fact{\Z[G]}{I_{\aug}^2}
$, and
$
I=\fact{I_{\aug}}{I_{\aug}^2}
$,
where $I_{\aug}$ is the augmentation ideal. Further, set
\begin{align*}
&\map{\pi}{\tilde{S}}{S}\\
&a_{\lambda,\mu}\mapsto \ovl{(g_{\gen})_{\lambda,\mu}}-\delta_{\lambda,\mu}\tp
\end{align*}

We prove that the pair $S$, $I$ satisfies the assumptions in the previous lemma. The
first two conditions are obvious, the third one is fulfilled because the augmentation homomorphism splits; it remains to show that $I$ is torsion free as an Abelian group. The
scheme $G(\Phi,-)$ is smooth, the homomorphism $\pi$ is surjective, and $(\Ker\pi)^2=0$; hence
the reduction homomorphism
$$
\map{\pi_*}{G(\Phi,\tilde{S})}{G(\Phi,S)}
$$
is surjective (see~\cite[\S4, Item 4.6]{DemazureGabriel}, set $k=\Z$). Hence there exists a matrix 
$$
e+(b_{\lambda,\mu})\in G(\Phi,\tilde{S})
$$
such that $\pi(\delta_{\lambda,\mu}+b_{\lambda,\mu})= (\ovl{g_{\gen}})_{\lambda,\mu}$. In particular, $b_{\lambda,\mu}\in \tilde{I}$; hence the ring homomorphism
\begin{align*}
\map{i&}{S}{\tilde{S}}\\
&(\ovl{g_{\gen}})_{\lambda,\mu}\mapsto b_{\lambda,\mu}+\delta_{\lambda,\mu}\tp
\end{align*}
is well defined. Then the homomorphism of Abelian groups
$$
\map{\restr{i}{I}}{I}{\tilde{I}}
$$
is a right inverse to $\restr{\pi}{\tilde{I}}$; hence $I$ is isomorphic to a direct summand of $\tilde{I}$ and, therefore, is tortion free

The assumption implies that the ring homomorphism
\begin{align*}
\map{\psi&}{S}{R}\\
&(\ovl{g_{\gen}})_{\lambda,\mu}\mapsto g_{\lambda,\mu}
\end{align*}
is well defined.

Now let weights $\lambda_0$,$\mu_0$ be such that $d(\lambda_0,\mu_0)\ge 2$. We prove that $g_{\lambda_0,\mu_0}=0$. To do this, it suffices to prove that $(\ovl{g_{\gen}})_{\lambda_0,\mu_0}=0$. Suppose the contrary; then by the previous
lemma, there exists a ring homomorphism
$$
\map{\ph}{S}{\fact{\C[\eps]}{(\eps^2)}}
$$
such that $\ph((\ovl{g_{\gen}})_{\lambda,\mu})\ne 0$.
The matrix $\ovl{g_{\gen}}$ belongs to $G(\Phi,S,I)$. Since $I^2=0$, we see that $\ph(I)$ belongs to the ideal generated by $\eps$; hence
$$
\ph_*(\ovl{g_{\gen}})\in G(\Phi,\fact{\C[\eps]}{(\eps^2)},(\eps))=\Lie(G(\Phi,-)_\C)\tp
$$ 
This Lie algebra admits a Cartan decomposition, i.e., the matrix $(\ph(a_{\lambda,\mu}))$ is a linear combination of elementary root elements and diagonal elements of this algebra. Thus,
since $d(\lambda_0,\mu_0)\ge 2$, its entry $\ph((\ovl{g_{\gen}})_{\lambda,\mu})$ is equal to 0, which contradicts the choice of~$\ph$. 
\end{proof}

 {\prop\label{NilpotentExtraction} 
  Let $H$ be an overgroup of $E(\Delta,R)$.
 \begin{enumerate}
 \item
 Let
$
(H\cap G(\Phi,R,\B)) \sm P^-\ne\emp\tc
$
where $\B\unlhd R$ is such that $\B^2=0$. Then $H$ contains an element of the form $x_\alpha(\xi)$, where $\alpha\in\Omega^+$ and $\xi\ne 0$.

\item The same holds for $P$ and $\Omega^-$. 
\end{enumerate} 
\begin{proof}
We prove Item (1); the proof of Item (2) is similar.

Let 
$
g\in (H\cap G(\Phi,R,\B)) \sm P^-\tp
$
Since $g\notin P^-$, there exists $\lambda_1\in\Lambda\sm\{\lambda_0\}$ such that
$g_{\lambda_0,\lambda_1}\ne 0$. By the previous lemma $\lambda_1\in\Lambda_1$. Take $\alpha\in\Delta$ such that $\lambda_1+\alpha=\nu\in\Lambda_1$(Lemma~\ref{weights}). If we prove that
$$
h=gx_\alpha(1)g^{-1}\in (H\cap P_{\lambda_1})\sm P^-\tc
$$
then by Proposition~\ref{PlambdaExtraction}, we are done.

We prove that the element $x_\alpha(1)$ stabilizes the line spanned by the vector $(g^{-1})_{*,\lambda_1}$.

Indeed, we have:
$$
(x_\alpha(1)(g^{-1})_{*,\lambda_1})_{\mu}=\begin{cases}
(g^{-1})_{\mu,\lambda_1}\quad\quad &\mu+\alpha\notin\Phi\tc\\
(g^{-1})_{\mu,\lambda_1}\pm (g^{-1})_{\mu+\alpha,\lambda_1}\quad\quad &\mu+\alpha\in\Phi\tp
\end{cases}
$$  
Note that since $\lambda_1+\alpha\in\Lambda$, in the second case, by Lemma~\ref{notroot} we have either $\mu=\lambda_1$ or $d(\mu+\alpha,\lambda_1)\ge 2$; hence $(g^{-1})_{\mu+\alpha,\lambda_1}=0$ by the previous lemma. Therefore,
$$
(x_\alpha(1)(g^{-1})_{*,\lambda_1})_{\mu}=\begin{cases}
(g^{-1})_{\mu,\lambda_1}\quad\quad &\mu\ne\lambda_1\tc\\
(g^{-1})_{\lambda_1,\lambda_1}\pm (g^{-1})_{\nu,\lambda_1}\quad\quad &\mu=\lambda_1\tp
\end{cases}
$$
Hence since $\B^2=0$, we obtain
$$
x_\alpha(1)(g^{-1})_{*,\lambda_1}=(g^{-1})_{*,\lambda_1}\pm (g^{-1})_{\nu,\lambda_1}v^{\lambda_1}=(g^{-1})_{*,\lambda_1}(1\pm (g^{-1})_{\nu,\lambda_1})\tp
$$
Thus $h\in P_{\lambda_1}$.

Finally, let $w$ be a row such that
$$
g_{\lambda_0,*}x_\alpha(1)=g_{\lambda_0,*}+w\tp
$$
Then $w_\mu\in\B$ for all $\mu\in\Lambda$ (note that $\lambda_0-\alpha\notin\Lambda$ because $\alpha\in\Delta$). Furthermore,
$
w_\nu=\pm g_{\lambda_0,\lambda_1}\ne 0\tp
$
Then since $g^{-1}\in G(\Phi,R,\B)$, we obtain:
$$
h_{\lambda_0,\nu}=(g_{\lambda_0,*}g^{-1}+wg^{-1})_\nu=0+w_\nu\ne 0\tp
$$
Thus 
$
h\notin P^-\tp
$
\end{proof}
}

\section{An $A_2$-proof}
The idea of the argument below is borrowed from \cite{VavGav}.

{\lem\label{A2} Let $g\in\Root(R)$ and $\lambda_1\in\Lambda$. For every pair of weights $\mu$,$\nu\in\Lambda$ such that $d(\lambda_1,\mu)=d(\lambda_1,\nu)=d(\mu,\nu)=1$, set
$$
x(\mu,\nu)=x_{\alpha}(c_{\mu,\alpha}\cdot g_{\nu,\lambda_1})x_{\beta}(-c_{\nu,\beta}\cdot g_{\mu,\lambda_1})\tc
$$
where $\alpha=\lambda_1-\mu$ and $\beta=\lambda_1-\nu$.

Then the following holds.
\begin{enumerate}
\item The element $x(\mu,\nu)$ stabilises the column $g_{*,\lambda_1}$.

\item We have
$$
g^{-1}x(\mu,\nu)g=x(\mu,\nu)+[x(\mu,\nu),g]_{\ring}\tc
$$
where $[\cdot,\cdot]_{\ring}$ is a ring commutator of matrices.
\end{enumerate}
 }
\begin{proof}
\begin{enumerate}
\item We act by the element $x_{\beta}(-c_{\nu,\beta}\cdot g_{\mu,\lambda_1})$ on the column $g_{*,\lambda_1}$
$$
(x_{\beta}(-c_{\nu,\beta}\cdot g_{\mu,\lambda_1})g_{*,\lambda_1})_\rho=\begin{cases}
g_{\rho,\lambda_1}\quad & \rho-\beta\notin\Lambda\tc\\
g_{\rho,\lambda_1}-c_{\rho-\beta,\beta}c_{\nu,\beta}\cdot g_{\mu,\lambda_1}g_{\rho-\beta,\lambda_1} \quad & \rho-\beta\in\Lambda\tp
\end{cases}
$$
Note that if the second case occurs, Lemma~\ref{notroot} shows that we have either $\rho=\lambda_1$ or $d(\rho-\beta,\lambda_1)\ge 2$; hence by Lemma~\ref{rootequation} $g_{\rho-\beta,\lambda_1}=0$. Therefore, we have:
$$
x_{\beta}(-c_{\nu,\beta}\cdot g_{\mu,\lambda_1})g_{*,\lambda_1}=g_{*,\lambda_1}-g_{\mu,\lambda_1}g_{\nu,\lambda_1}v^{\lambda_1}\tp
$$
Similarly, applying $x_{\alpha}(c_{\mu,\alpha}\cdot g_{\nu,\lambda_1})$ to the resulting column, we obtain the initial column $g_{*,\lambda_1}$.

\item It suffices to consider the case where $R=\Z[\Root]$ and $g=r_{\gen}$. Since what we need
to prove is a polynomial identity, we may extend the ring up to $\Z[G]$, using Lemma~\ref{ZRootToZG}, and then extend it up to $\Z[G][{1\over 2}]$.

Let $\tilde{x}=x(\mu,\nu)-e$ and $\tilde{r}=r_{\gen}-e$ be the corresponding elements of the Lie algebra (see Lemma~\ref{RootTypeAsExponent}).  Here it should be noted that $\angle(\alpha,\beta)=\tfrac{\pi}{3}$ because 
$
\alpha-\beta=\nu-\mu\in\Phi\tp
$
Hence $x(\mu,\nu)\in\Root(R)$ by Lemma~\ref{examples}. Then by Lemma~\ref{RootTypeAsExponent}, we have the following identity:
\begin{align*}
r_{\gen}^{-1}x(\mu,\nu)r_{\gen}=e+\exp(-\ad \tilde{r})\tilde{x}=x(\mu,\nu)+[\tilde{x},\tilde{r}]_{\ring}+{1\over2}(\ad\tilde{r})^2\tilde{x}=x(\mu,\nu)+\\+[x(\mu,\nu),r_{\gen}]_{\ring}+{1\over2}(\ad\tilde{r})^2\tilde{x}\tp
\end{align*}
Next, note that for any element $y$ of our Lie algebra, we have
$
(\ad e_{\alpha_1})^2y=\xi_y e_{\alpha_1}
$
for some $\xi_y\in R$, where $e_{\alpha_1}=x_{\alpha_1}(1)-e$  is a root element. This relation is obtained from the Cartan decomposition of the element $y$ by expanding (the only nonzero term comes from the coefficient of $e_{-\alpha_1}$). Conjugating this relation by the element $g_{\gen}$,  we obtain a similar statement for $r_{\gen}$. Therefore, we have the formula
$$
r_{\gen}^{-1}x(\mu,\nu)r_{\gen}=x(\mu,\nu)+[x(\mu,\nu),r_{\gen}]_{\ring}+\xi\tilde{r}\eqno(*)
$$
for some $\xi\in \Z[G][{1\over 2}]$. We claim that in fact, we have $\xi=0$. Indeed, let us compare the $\lambda_1$th columns of the left-hand and right-hand sides of $(*)$. The first item implies the formula
$$
(r_{\gen}^{-1}x(\mu,\nu)r_{\gen})_{*,\lambda_1}=v_{\lambda_1}\tc
$$ 
as well as the relation
$$
(x(\mu,\nu)r_{\gen})_{*,\lambda_1}=(r_{\gen})_{*,\lambda_1}\tp
$$ 
Further, it is easy to see that the element $x(\mu,\nu)$ also stabilises the column $v_{\lambda_1}$. Hence
$
x(\mu,\nu)_{*,\lambda_1}=v_{\lambda_1}\tp
$ 
Finally, calculating the product $r_{gen}x(\mu,\nu)$ we see that the matrix $x(\mu,\nu)$ acts on each row of the matrix $r_{gen}$, without changing the entry in the $\lambda_1$th column (indeed, $\lambda_1-\alpha$ and $\lambda_1-\beta$ belong to $\Lambda$; hence  $\lambda_1+\alpha$ and $\lambda_1+\beta$ do not belong to $\Lambda$,  because the representation is minuscule), i.e., we have
$$
(r_{gen}x(\mu,\nu))_{*,\lambda_1}=(r_{\gen})_{*,\lambda_1}\tp
$$
Therefore, we obtain the identity:
$$
v_{\lambda_1}=v_{\lambda_1}+(r_{\gen})_{*,\lambda_1}-(r_{\gen})_{*,\lambda_1}+\xi\tilde{r}_{*,\lambda_1}\tc
$$
i.e.,
$
\xi\tilde{r}_{*,\lambda_1}=0\tp
$

The ring $\Z[G][{1\over 2}]$ has no zero divisors, and the column $\tilde{r}_{*,\lambda_1}$ is nonzero. Hence $\xi=0$, and $(*)$ gives us what we need.
\end{enumerate}
\end{proof}

{\lem\label{A1} Let $g\in\Root(R)$ and $\lambda_1\in\Lambda$. Let $\mu\in\Lambda$ be such that $d(\lambda_1,\mu)=1$, and let $\xi\in R$  be such that $\xi g_{\mu,\lambda_1}=0$. Then the element $x_\alpha(\xi)$, where $\alpha=\lambda_1-\mu$, stabilises the column $g_{*,\lambda_1}$.}
\begin{proof} We have
$$
(x_{\alpha}(\xi)g_{*,\lambda_1})_\rho=\begin{cases}
g_{\rho,\lambda_1}\quad & \rho-\alpha\notin\Lambda\tc\\
g_{\rho,\lambda_1}\pm \xi g_{\rho-\alpha,\lambda_1} \quad & \rho-\alpha\in\Lambda\tc
\end{cases}
$$
moreover, if the second case occurs, Lemma~\ref{notroot} shows that either $\rho=\lambda_1$, or $d(\rho-\alpha,\lambda_1)\ge 2$; hence $g_{\rho-\alpha,\lambda_1}=0$ by Lemma~\ref{rootequation}. Therefore, we obtain
$$
x_{\alpha}(\xi)g_{*,\lambda_1}=g_{*,\lambda_1}\pm \xi g_{\mu,\lambda_1}v^{\lambda_1}=g_{*,\lambda_1}\tp
$$
\end{proof}

We introduce some notation. For each element $g\in \Root(R)$ and weight $\lambda_1\in \Lambda_1$ consider the following four ideals of the ring $R$:
\begin{align*}
&\A(g,\lambda_1)=\<\{g_{\mu,\lambda_1}\colon \mu\in\Lambda_1\sm\{\lambda_1\}\}\>\tc\\
&\B(g,\lambda_1)=\<g_{\lambda_0,\lambda_1}\>\tc\\
&\A'(g,\lambda_1)=\<\{g_{\lambda_1,\mu}\colon \mu\in\Lambda_1\sm\{\lambda_1\}\}\>\tc\\
&\B'(g,\lambda_1)=\<g_{\lambda_1,\lambda_0}\>\tp\\
\end{align*}

{\prop\label{A2proof} Let $H$ be an overgroup of $E(\Delta,R)$ of level $\sigma=(I^+,I^-)$, let $g\in\Root(R)\cap H$, and let $\lambda_1\in\Lambda_1$. Then we have:
\begin{enumerate}
\item $\A(g,\lambda_1)\B(g,\lambda_1)\sub I^+$;
\item if $I^+=0$, then $(\B(g,\lambda_1))^3=(0)$;
\item $\A'(g,\lambda_1)\B'(g,\lambda_1)\sub I^-$;
\item if $I^-=0$, then $(\B'(g,\lambda_1))^3=(0)$.
\end{enumerate}}
\begin{proof}
\begin{enumerate}
\item Let $\nu\in \Lambda_1\sm \{\lambda_1\}$, We prove that $g_{\nu,\lambda_1}g_{\lambda_0,\lambda_1}\in I^+$. It suffices to consider the case where $d(\nu,\lambda_1)=1$ because otherwise $g_{\nu,\lambda_1}=0$ by Lemma~\ref{rootequation}. 

Take a weight $\mu$ from Item 2 of Lemma~\ref{weights}. Let $\alpha=\lambda_1-\mu$, $\beta=\lambda_1-\nu$ and let $x(\mu,\nu)$ be as in Lemma~\ref{A2}. Then this lemma implies that the element
$
g^{-1}x(\mu,\nu)g
$
stabilises $v^{\lambda_1}$ (i.e., belongs to $P_{\lambda_1}$).

By Lemma~\ref{examples} this element belongs to $\Root(R)$. It also belongs to $H$ (indeed,
since $\lambda_1$, $\mu$, $\nu\in\Lambda_1$, we have $\alpha$, $\beta\in\Delta$; hence $x(\mu,\nu)\in H$). If it does not belong to $G(\Phi,\Delta,R,(I^+,R))$, then by Proposition~\ref{PlambdaExtraction} the group $H$ cannot be of level $\sigma$. Hence it belongs to $G(\Phi,\Delta,R,(I^+,R))$, and by Item 2 of Lemma~\ref{A2} we obtain the congruence
$$
(gx(\mu,\nu))_{\lambda_0,\mu}\equiv (x(\mu,\nu)g)_{\lambda_0,\mu}\mod I^+\tp\eqno{(1)}
$$
We compute the left-hand side of the congruence (1). Note that $\mu-\beta\notin\Lambda$, because otherwise, by Lemma~\ref{notroot}, the relation $\lambda_1+\beta\in\Lambda$, implies that either $d(\lambda_1,\mu)\ge 2$, or $\mu=\nu$, which contradicts the choice of $\mu$. Hence,
$$
gx(\mu,\nu)v^\mu=gx_{\alpha}(\pm g_{\nu,\lambda_1})x_{\beta}(\pm g_{\mu,\lambda_1})v^\mu=gx_{\alpha}(\pm g_{\nu,\lambda_1})v^\mu=g(v^\mu\pm g_{\nu,\lambda_1}v^{\lambda_1})\tc
$$
which implies that
$$
(gx(\mu,\nu))_{\lambda_0,\mu}=g_{\lambda_0,\mu}\pm g_{\nu,\lambda_1}g_{\lambda_0,\lambda_1}\tp
$$
Now we compute the right-hand side of the congruence (1). Note that $\lambda_0-\alpha\notin\Lambda$ because $\alpha\in\Delta$.  Hence,
$$
(x_{\alpha}(\pm g_{\nu,\lambda_1})x_{\beta}(\pm g_{\mu,\lambda_1})g)_{\lambda_0,\mu}=(x_{\beta}(\pm g_{\mu,\lambda_1})g)_{\lambda_0,\mu}\tp
$$
Similarly, $\lambda_0-\beta$ does not belong to $\Lambda$; hence we obtain
$$
(x_{\beta}(\pm g_{\mu,\lambda_1})g)_{\lambda_0,\mu}=g_{\lambda_0,\mu}\tp
$$
Therefore, the congruence (1) tells us exactly that $g_{\nu,\lambda_1}g_{\lambda_0,\lambda_1}\in I^+$.

\item We prove that $(g_{\lambda_0,\lambda_1})^3=0$. 

Take a weight $\mu$ from Item 1 of Lemma~\ref{weights}. Put $\alpha=\lambda_1-\mu$, and $x=x_\alpha(g_{\lambda_0,\lambda_1})$. By what we proved before, we have $g_{\lambda_0,\lambda_1}g_{\mu,\lambda_1}=0$. Then Lemma~\ref{A1} implies
that the element
$
g^{-1}xg
$
stabilizes $v^{\lambda_1}$.  It also belongs to $H\cap\Root(R)$. If it does not
belong to $P^-$ (i.e., to $G(\Phi,\Delta,R,(0,R))$), , then by Proposition~\ref{PlambdaExtraction} the group $H$ cannot be of level $\sigma$ (because $I^+=0$). Hence it belongs to $P^-$, i.e., the matrix  $x$ stabilizes the line spanned by the covector $g^{-1}_{\lambda_0,*}$.  Thus it multiplies this covector by a scalar, which we denote by $(1+\eps)$. Then, first, we have:
$$
(1+\eps)(g^{-1})_{\lambda_0,\lambda_1}=(g^{-1})_{\lambda_0,\lambda_1}\tc
$$
i.e., 
$
\eps (g^{-1})_{\lambda_0,\lambda_1}=0\tp
$
Second, we have:
$$
(1+\eps)(g^{-1})_{\lambda_0,\mu}=(g^{-1}x)_{\lambda_0,\mu}=(g^{-1})_{\lambda_0,\mu}\pm (g^{-1})_{\lambda_0,\lambda_1}^2\tc
$$
i.e.,
$
\eps(g^{-1})_{\lambda_0,\mu}=\pm (g^{-1})_{\lambda_0,\lambda_1}^2\tp
$
Therefore, we obtain the identity
$$
(g^{-1})_{\lambda_0,\lambda_1}^3=\pm\eps(g^{-1})_{\lambda_0,\mu}(g^{-1})_{\lambda_0,\lambda_1}=0\tp
$$
It remains to note that by Lemma~\ref{RootTypeAsExponent}, we have:
$
(g^{-1})_{\lambda_0,\lambda_1}=-g_{\lambda_0,\lambda_1}\tp
$
\end{enumerate}

The proofs of Items 3 and 4 are similar. 
\end{proof}

{\cor\label{lift} Let $I$ be an ideal of $R$, and let $H$ be an overgroup of $E(\Delta,R)$ of level $\sigma=(I^+,I^-)$. Then for the level of $\rho_I(H)$ as an overgroup of $E\left(\Delta,\fact{R}{I}\right)$ we have:
$$
\lev(\rho_I(H))=(\rho_I(I^+),\rho_I(I^-))\tp
$$}
\begin{proof}
Obviously, the right-hand side is contained in the left-hand side, let us prove the
inverse inclusion.

Suppose the contrary. Without loss of generality, we may assume that the inclusion
fails for the first component, i.e., for some $g\in H$, $\ovl{\xi}\in \fact{R}{I}$ such that $\ovl{\xi}\notin \rho_I(I^+)$, and $\alpha\in\Omega^+$ we have
$
\rho_I(g)=x_\alpha(\ovl{\xi})\tp
$

Let $\nu=\lambda_0-\alpha\in\Lambda_1$ (Lemma~\ref{OmegaAndHighestWeight}). Take $\beta\in\Delta$ such that $\alpha+\beta\in\Phi$ (Lemma~\ref{RootInDeltaComb}). Then $\alpha+\beta\in\Omega^+$, and $\lambda_1=\lambda_0-\alpha-\beta\in\Lambda_1$.

Set 
$
h=(x_\beta(1))^g\tp
$
Clearly, we have $h\in H\cap \Root(R)$. Next, we have
$$
\rho_\B(h)=(x_\beta(1))^{x_\alpha(\ovl{\xi})}=x_\beta(-1)\cdot[x_\beta(1),x_\alpha(-\ovl{\xi})]=x_\beta(-1)x_{\alpha+\beta}(\pm\ovl{\xi})
$$
By the previous proposition, we have
$
h_{\nu,\lambda_1}h_{\lambda_0,\lambda_1}\in I^+\tp
$
Hence we have
$$
\rho_\B(h)_{\nu,\lambda_1}\rho_\B(h)_{\lambda_0,\lambda_1}\in\rho_I(I^+)\tp
$$
On the other hand, the following identity holds:
\begin{align*}
\rho_\B(h)v^\lambda=x_\beta(-1)x_{\alpha+\beta}(\pm\ovl{\xi})v^\lambda=x_\beta(-1)(v^\lambda\pm \ovl{\xi}v^{\lambda+\alpha+\beta})=\\=\begin{cases}
v^\lambda\pm \ovl{\xi}v^{\lambda+\alpha+\beta}\pm v^{\lambda+\beta} & \lambda+\alpha+2\beta\notin\Lambda\tc\\
v^\lambda\pm \ovl{\xi}v^{\lambda+\alpha+\beta}\pm v^{\lambda+\beta}\pm \ovl{\xi}v^{\lambda+\alpha+2\beta} & \lambda+\alpha+2\beta\in\Lambda\tp
\end{cases}
\end{align*}

In any case, we obtain $\rho_\B(h)_{\lambda+\beta+\alpha,\lambda}=\pm\ovl{\xi}$, and $\rho_\B(h)_{\lambda+\beta,\lambda}=\pm 1$. Hence $\ovl{\xi}\in\rho_I(I^+)$, which contradicts the assumption.
\end{proof}

\section{Root type elements in an overgroup}
{\lem\label{RootInDelta} Let $\sigma=(I^+,I^-)$ be a pair of ideals, and let
 $g\in \Root(R)\sm G(\Phi,\Delta,R,\sigma)$. Then there exists $\gamma\in \Delta$ such that 
$$
gx_{\gamma}(1)g^{-1}\notin G(\Phi,\Delta,R,\sigma)\tp
$$}
\begin{proof}
Clearly, without loss of generality we may assume that 
$
g\notin G(\Phi,\Delta,R,(I^+,R))\tp
$
 First we prove the following statement:
 $$
 gE(\Delta,R)g^{-1}\not\sub G(\Phi,\Delta,R,(I^+,R))\tp
 $$
 
 Suppose the contrary; then our assumption means exactly that the group $E(\Delta,R)$ stabilizes the line spanned by the covector $g_{\lambda_0*}$ modulo $I^+$.  Since this group is perfect, it must stabilize the covector itself.
  
Further since $g\notin G(\Phi,\Delta,R,(I^+,R))$, there exists $\mu\in\Lambda\sm\{\lambda_0\}$ such that $g_{\lambda_0,\mu}\notin I^+$. By Lemma~\ref{rootequation}, we have $\mu\in\Lambda_1$. Take $\tilde{\gamma}\in\Delta$ such that
  $$
  \nu=\mu+\tilde{\gamma}\in \Lambda_1
  $$
  (apply Lemma \ref{RootInDeltaComb} for $\beta=\lambda_0-\nu$).
  
  Since the element  $x_{\tilde{\gamma}}(1)$ stabilizes the covector $g_{\lambda_0*}$, we obtain the following congruence:
  $$
  (gx_{\tilde{\gamma}}(1))_{\lambda_0,\nu}\equiv g_{\lambda_0,\nu}\mod I^+\tp
  $$
 However, we also have
  $$
  (gx_{\tilde{\gamma}}(1))_{\lambda_0,\nu}=g_{\lambda_0,\nu}\pm g_{\lambda_0,\mu}\tp
  $$
  This implies that $g_{\lambda_0,\mu}\in I^+$, which contradicts the choice of $\mu$.
  
  Therefore, we have proved that for some $\gamma\in\Delta$ and some $t\in R$,  we have
  $$
 gx_{\gamma}(t)g^{-1}\notin G(\Phi,\Delta,R,(I^+,R))\tp
 $$
 We also know that
 $$
 gx_{\gamma}(t)g^{-1}=e+t(gx_{\gamma}(1)g^{-1}-e)\tc
 $$
which implies
 $$
 gx_{\gamma}(1)g^{-1}\notin G(\Phi,\Delta,R,(I^+,R))\tp
 $$
\end{proof}  

{\prop\label{RootElementsInH} Let $H$ be an overgroup of $E(\Delta,R)$ of level $\sigma=(I^+,I^-)$, and let $g\in H\cap \Root(R)$. Then $g\in G(\Phi,\Delta,R,\sigma)$.}
\begin{proof}
First, we reduce the proof to the case where one of the ideals in $\sigma$ is equal to zero.
Assume that we can prove the proposition for such $\sigma$, and let us prove it in general case.
Note that
$
\rho_{I^+}(g)\in \rho_{I^+}(H)\cap \Root(\fact{R}{I^+})\tp
$
By Corollary~\ref{lift} we have:
$$
\lev(\rho_{I^+}(H))=((0),\rho_{I^+}(I^-))\tp
$$
Hence by assumption
$$
\rho_{I^+}(g)\in G(\Phi,\Delta,\fact{R}{I^+},\ldots)\le P^-\tp
$$
Similarly, we obtain
$
\rho_{I^-}(g)\in P\tc
$
q.e.d.

Therefore, without loss of generality, we may assume that $I^+=0$. Now we prove the proposition in this case.

Suppose the contrary. Take $\alpha$ from Lemma~\ref{RootInDelta}, applied to the element $g^{-1}$, i.e.,
$$
g^{-1}(x_\alpha(1))g\notin G(\Phi,\Delta,R,\sigma)\tp
$$
Set
$$
M=\{\lambda_1\in\Lambda_1\colon g_{\lambda_0,\lambda_1}\ne 0\}\tp
$$
For any $\lambda_1\in M$, set
$$
k(\lambda_1)=\min\{k\colon (g_{\lambda_0,\lambda_1})^k =0\}
$$
(by Proposition~\ref{A2proof}, we have $k(\lambda_1)\le 3$). Finally, set
$$
K=\sum_{\lambda_1\in M} k(\lambda_1)\tp
$$

We prove by induction on $K$ that our assumption leads to a contradiction.

The base of induction: $K=0$; this means that $M=\emp$. Together with Item 1 of
Lemma~\ref{rootequation} this means exactly that $g\in P^-$.
Then $g\in G(\Phi,\Delta,R,\sigma)$ because otherwise,
by Item 2 of Proposition~\ref{PExtraction}, the group $H$ cannot be of level $\sigma$.

Now we pass to the induction step. Assume that $M\ne\emp$. Fix $\lambda_1\in M$. Set 
$$
\B=\begin{cases}
\B(g,\lambda_1)\quad k(\lambda_1)=2\\
\B(g,\lambda_1)^2\quad k(\lambda_1)=3
\end{cases}
$$
Then $\B^2=(0)$.

{\bf Case $1$}. $\rho_\B(g)\notin G(\Phi,\Delta,\fact{R}{\B},\rho_\B(\sigma))$.

Note that 
$
\rho_{\B}(g)\in \rho_{\B}(H)\cap \Root(\fact{R}{\B})\tp
$
By Corollary~\ref{lift}
$
\lev(\rho_{\B}(H))=\rho_\B(\sigma)\tc
$
next, the number $K$ for $\rho_\B(g)$ is less than that for $g$. Applying the inductive hypothesis, we obtain a contradiction.

{\bf Case $2$}. $\rho_\B(g)\in G(\Phi,\Delta,\fact{R}{\B},\rho_\B(\sigma))$.

By Lemma~\ref{normality},
$$
(x_\alpha(1))^{\rho_\B(g)}\in E(\Phi,\Delta,\fact{R}{\B},\rho_\B(\sigma))\tp
$$
The group $ E(\Phi,\Delta,R,\sigma)$ ) maps surjectively onto the group $E(\Phi,\Delta,\fact{R}{\B},\rho_\B(\sigma))$. Hence there exists $h\in E(\Phi,\Delta,R,\sigma)$ such that
$$
\rho_\B(h)=(x_\alpha(1))^{\rho_\B(g)}=\rho_\B((x_\alpha(1))^g)\tp
$$
Then we have
$$
g_1=h^{-1}(x_\alpha(1))^g\in (H\cap G(\Phi,R,\B))\sm G(\Phi,\Delta,R,\sigma)\tp
$$
If $g_1\in P^-$, then the group $H$ cannot be of level $\sigma$ by Proposition~\ref{PExtraction}, and otherwise it
cannot by Proposition~\ref{NilpotentExtraction}. This is a contradiction. 
\end{proof}

\section{Finishing the proof of Theorem~1}
Let $H$ be an overgroup of $E(\Delta,R)$ of level $\sigma=(I^+,I^-)$. We want to prove that
$$
H\le N_{G(\Phi,R)}\left(E(\Phi,\Delta,R,\sigma)\right)\tp
$$
By Proposition~\ref{Transporter}, it suffices to prove that for any $g\in H$ and any generator $x_\alpha(\xi)$ of the group $E(\Phi,\Delta,R,\sigma)$ we have
$
g^{-1}x_\alpha(\xi)g\in G(\Phi,\Delta,R,\sigma)\tp
$
which is true indeed by Proposition~\ref{RootElementsInH}.

\section{Changing the lattice~$\P$}

In this section we generalise Theorem~\ref{sandwich} to other weight lattices. First note that Theorem~\ref{sandwich} and Corollary~\ref{DerivSer2} imply the following fact.

{\cor\label{DerivSer3} Let $R$ R be a finitely generated ring, and let $H$ be an overgroup of $E(\Delta,R)$ of level $\sigma$. Then there exists a natural number $N$ such that
$$
D^NH=E(\Phi,\Delta,R,\sigma)\tp
$$}

Note also that Corollary~\ref{DerivSer3} is a stronger statement than Theorem~\ref{sandwich}. Indeed, for finitely generated rings, Theorem~\ref{sandwich}  follows from Corollary~\ref{DerivSer3} directly, and the general case can
be derived by using the fact that any ring is an inductive limit of its finitely generated
subrings. Let us transfer the result of Theorem~\ref{sandwich} to other weight lattices in this stronger form

Let $\P'$ be another weight lattice between $\euQ(\Phi)$ and $\P(\Phi)$. We denote by $G_{\P'}(\Phi,R)$ the corresponding Chevalley group. The subgroups $E_{\P'}(\Delta,R)$ and $E_{\P'}(\Phi,\Delta,R,\sigma)$ for it are defined similarly. In the case of $\P'=\P(\Phi)$ we will use the subscript $\SC$. Let $C_{\SC}(\Phi,R)$  be the center of the group $G_{\SC}(\Phi,R)$.

For any intermediate subgroup
$$
E_{\P'}(\Delta,R)\le H\le G_{\P'}(\Phi,R)
$$
one can similarly define its level.

{\prop Corollary~$\ref{DerivSer3}$, and hence Theorem~$\ref{sandwich}$, can be transferred to the simply connected group.
}
\begin{proof}
Let $H$ be an overgroup of $E_{\SC}(\Phi,R)$ of level $\sigma=(I^+,I^-)$. 
 
Consider the natural homomorphism
$$
\map{\pi}{G_{\SC}(\Phi,R)}{G(\Phi,R)}\tp
$$
Obviously, we have $\lev\pi(H)\ge\sigma$.  Let us prove the inverse inclusion. Assume that
$
x_\beta(\xi)\in\pi(H)\tp
$
for some $\beta\in\Omega^+$ and $\xi\in R\sm I^+$ (the proof for $\Omega^-$ is similar). Then for some $g\in C_{\SC}(\Phi,R)$ we have
$
gx_\beta(\xi)\in H\tp
$
Take $\alpha$ from Lemma~\ref{RootInDeltaComb}, then we have
$$
x_{\alpha+\beta}(\pm\xi)=[x_\alpha(1),x_{\beta}(\xi)]=[x_\alpha(1),gx_{\beta}(\xi)]\in H\tc
$$
which contradicts the assumption

Thus $\lev\pi(H)=\sigma$. Hence for a large $N$ the following inclusion holds
$$
D^N\pi(H)\le E(\Phi,\Delta,R,\sigma)\tp
$$
Therefore, we have the inclusion
$$
D^N H\le C_{\sc}(\Phi,R)E_{\sc}(\Phi,\Delta,R,\sigma)\tc
$$
which, in its turn, implies that
$$
D^{N+1} H\le E_{\sc}(\Phi,\Delta,R,\sigma)\tp
$$
\end{proof} 

{\prop \label{Pprime}
  Corollary $\ref{DerivSer3}$, and hence Theorem~$\ref{sandwich}$, can be transfered to an arbitrary weight lattice $\P'$.
}
\begin{proof}
Let
$$
\map{\pi}{G_{\SC}}{G_{\P'}(\Phi,R)}
$$
be the natural homomorphism.

The sequence
$$
\xymatrix{
1\ar@{->}[r] & \Ker\pi\ar@{->}[r] & G_{\SC}(\Phi,-)\ar@{->}[r]^{\pi} &  G_{\P'}(\Phi,-)\ar@{->}[r] & 1\tp
}
$$
is exact in the category of fppf-sheaves. Hence for groups we have the exact sequence
$$
\xymatrix{
1\ar@{->}[r] & \Ker\pi\ar@{->}[r] & G_{\SC}(\Phi,R)\ar@{->}[r]^{\pi} &  G_{\P'}(\Phi,R)\ar@{->}[r] & H^1_{\fppf}(R,\Ker\pi)\tp
}
$$
\end{proof}

Therefore, $\Im\pi$ is a normal subgroup of $ G_{\P'}(\Phi,R)$ with an Abelian quotient (which
can be embedded into the Abelian group $ H^1_{\fppf}(R,\Ker\pi)$).

Now let $H$ be an overgroup of $E_{\P'}(\Delta,R)$  of level $\sigma$. Replacing $H$ with its commutant, we may assume that $H\sub \Im\pi$. Clearly, we have  $\lev(\pi^{-1}(H))=\sigma$. Then for a large $N$ we have
$$
D^N\pi^{-1}(H)=E_{\SC}(\Phi,\Delta,R,\sigma)\tp
$$
Since $H\le \Im\pi$, we have
$
H=\pi(\pi^{-1}(H))\tp
$
which implies the formula
$$
D^N H=\pi(D^N \pi^{-1}(H))=\pi(E_{\SC}(\Phi,\Delta,R,\sigma))=E_{\P'}(\Phi,\Delta,R,\sigma)\tp
$$

\end{document}